\documentclass[11pt]{article}
\usepackage[utf8]{inputenc}

\usepackage{bm}
\usepackage{amsmath}
\usepackage{amsfonts}
\usepackage{amssymb}
\usepackage{enumitem}
\usepackage{float}
\usepackage[subrefformat=parens,labelformat=parens,justification=centering]{subcaption}
\usepackage{upgreek}
\usepackage[numbers,sort&compress]{natbib}
\usepackage{xcolor}
\usepackage[margin=25mm,top=25mm]{geometry}
\usepackage{graphicx}
\usepackage{hyperref}

\definecolor{gray1}{rgb}{0.5,0.5,0.5}
\definecolor{col1}{rgb}{0,0.2353,0.6824}
\definecolor{col2}{rgb}{0.7804,0.0039,0.0549}
\definecolor{col3}{rgb}{0.2353,0.6078,0.2745}
\definecolor{col4}{rgb}{0.6,0.3,0.6}
\definecolor{col5}{rgb}{1,0.6078,0}

\graphicspath{{./figures/}}

\providecommand{\abstractn}[1]
{
  \begin{center}\textbf{Abstract}\\
  \vspace{5pt}
  \begin{minipage}{0.8\textwidth}
    \hspace{-\parindent}#1 
  \end{minipage}
  \end{center}
}

\providecommand{\keywords}[1]
{
    \begin{center}
  \begin{minipage}{0.8\textwidth}
    \small	
    \hspace{-\parindent}\textbf{Keywords: } #1 
  \end{minipage}
  \vspace{10pt}
  \end{center}
}

\title{\textbf{Modal-based synthesis of passive electrical networks for multimodal piezoelectric damping}}
\author{G. Raze$^{1}$, J. Dietrich$^{1}$, B. Lossouarn$^{2}$ and G. Kerschen$^{1}$\\
        \small $^{1}$Department of Aerospace and Mechanical Engineering, \\
        \small University of Liège, 4000 Liège, Belgium\\
        \small $^{2}$Laboratoire de Mécanique des Structures et des Systèmes Couplés,\\
        \small Conservatoire national des arts et métiers, 75003 Paris, France
        }
\date{}

\begin{document}

\maketitle

\abstractn{This work presents a new approach to design an electrical network which, when coupled to a structure through an array of piezoelectric transducers, provides multimodal vibration mitigation. The characteristics of the network are specified in terms of modal properties. On the one hand, the electrical resonance frequencies are chosen to be close to those of the targeted set of structural modes. On the other hand, the electrical mode shapes are selected to maximize the electromechanical coupling between the mechanical and electrical modes while guaranteeing the passivity of the network. The effectiveness of this modal-based synthesis is demonstrated using a free-free beam and a fully clamped plate.}

\keywords{vibration mitigation, piezoelectric damping, multimodal vibration absorber, passive control, electrical network, electromechanical analogy}

In several engineering applications, flexible and lightly-damped structures are more and more common but are plagued by harmful high-amplitude vibrations. Piezoelectric shunt damping such as proposed by Hagood and von Flotow~\cite{Hagood1991} and Wu~\cite{Wu1996} constitutes an attractive solution to this problem because of its compact and potentially lightweight character. In this approach, one or several piezoelectric transducers bonded to a structure convert part of their mechanical energy to electrical energy via the piezoelectric effect. This energy is then dissipated by electrical elements in a passive circuit connecting the electrodes of the transducer.

The approaches in~\cite{Hagood1991,Wu1996} can be used to mitigate a single resonant mode. This may become a limitation if the excitation is broadband and if multiple modes respond to it. Approaches that provide multimodal damping using a single piezoelectric transducer connected to a circuit with multiple branches were proposed, see, e.g.~\cite{Moheimani2006,Berardengo2017,Raze2020} for an overview. While using a single transducer makes up for a potentially compact solution, its placement on the structure may limit the performance in terms of vibration reduction of some modes. From this standpoint, it is desirable to use multiple transducers distributed over the structure to ensure sufficient electromechanical coupling with all the targeted modes. However, finding an electrical circuit that interconnects multiple piezoelectric transducers and optimally damps multiple modes is a challenging task.

Dell'Isola and Vidoli~\cite{DellIsola1998} used piezoelectric actuators uniformly distributed over a truss beam and designed an electrical transmission line with similar wavespeed as a target mechanical wave speed, opening the possibility for efficient electrical dissipation of the mechanical energy. Vidoli and Dell'Isola~\cite{Vidoli2000} then studied the spectral properties of the operators governing the dynamics of a beam and characterized the coupling existing between mechanical and electrical modes. An important conclusion of their work was that the mechanical and electrical modal characteristics should coincide in order to obtain strong coupling. This lead Alessandroni et al~\cite{Alessandroni2002} to find analog electrical circuits to beams and plates, i.e., having identical resonance frequencies and mode shapes. They also derived optimal resistances to efficiently damp one targeted mode. Maurini et al~\cite{Maurini2004} studied different circuit topologies and showed that one of them could provide multimodal damping of a beam. Porfiri et al~\cite{Porfiri2004} synthesized a passive circuit analog to a vibrating beam made only of inductors, capacitors and ideal transformers. The theoretical concepts presented in these works were experimentally validated on rods~\cite{Lossouarn2015}, beams~\cite{Lossouarn2015a} and plates~\cite{Lossouarn2016,Darleux2020a}. Reviews on this approach can be found in~\cite{Giorgio2015,Darleux2020}.

The aforementioned works used a homogenized model of a simple structure, later discretized with finite differences in order to derive a circuit with lumped elements being analog to the structure. Real-life structures may be more complex. This is what motivated Darleux et al~\cite{Darleux2020a,Darleux2020} to develop a library of elemental electrical cells analog to various elemental structures. By assembling these cells into a network in a similar way to that in which their mechanical counterparts make up a complex structure, it is possible to create a network analog to that structure. This approach potentially requires a large number of cells and transducers.

In a slightly different spirit, Giorgio et al~\cite{Giorgio2009} proposed a generic approach, where an electrical network was tuned based on the finite element model of a structure. The idea developed therein was to find a transformation of the electrical degrees of freedom of the network that makes the piezoelectric coupling matrix nearly diagonal in order to consider mechanical and electrical modes by pairs, thereby allowing a tuning based on the classical resistive-inductive shunt formulas~\cite{Wu1996}. This method was numerically verified and experimentally validated. However, it has two potential drawbacks. The first one is that it requires to solve a quadratic system of $N_s^2$ equations, where $N_s$ is the number of modes to be controlled. For such systems, an iterative numerical solver is needed which may not always converge to a satisfactory solution. The second issue is that the number of piezoelectric transducers needs to be equal to the number of targeted mechanical modes, which somewhat limits the flexibility of this approach.

The goal of the present paper is to propose a procedure for synthesizing an electrical network to be used for multimodal piezoelectric damping based on the modal properties of the host structure. The proposed approach specifies the characteristics of the resonant electrical modes of the network while requiring that the network is realizable using passive electrical elements. Specifically, by matching the electrical resonance frequencies to those of a set of targeted mechanical modes, and by analytically optimizing the electrical mode shapes, it is possible to derive the nodal admittance matrix of the network allowing to efficiently mitigate the targeted mechanical resonances. The proposed method, termed \textit{modal-based synthesis}, is non-iterative and can accommodate a different number of transducers and targeted modes.

This article is organized as follows. The basics of piezoelectric shunt damping are first briefly reviewed in \autoref{sec:piezoShunt}. An electrical network aiming to provide multimodal vibration damping is introduced in \autoref{sec:dynamics}, and the dynamics of the coupled system are analyzed. The electromechanical interaction between a mechanical mode and its electrical counterpart is described in terms of coupling factors. The actual design of the electrical network is undertaken in \autoref{sec:networkDesign}. Limitations on the electrical modal characteristics are imposed by the passivity requirement on the network. Under these constraints, the electrical modal characteristics are optimized and the network is synthesized. The approach is numerically demonstrated on a free-free beam and a fully clamped plate in \autoref{sec:examples}. Conclusions on the present work are finally drawn in \autoref{sec:conclusion}.

\section{Piezoelectric shunt damping}
\label{sec:piezoShunt}
        
        \begin{figure}[!ht]
			\centering
			\begin{subfigure}[c]{0.45\textwidth}
		         \centering
		         \includegraphics[width=\textwidth]{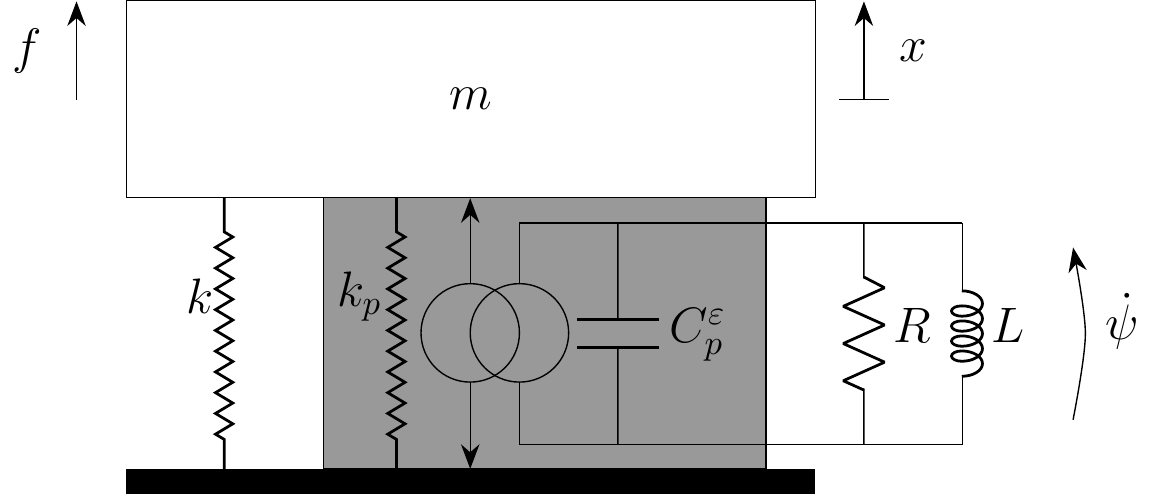}
		         \caption{}
		         \label{sfig:sdof_schematics}
		     \end{subfigure}
		     \hspace{0.05\textwidth}
			\begin{subfigure}[c]{0.4\textwidth}
		         \centering
		         \includegraphics[width=\textwidth]{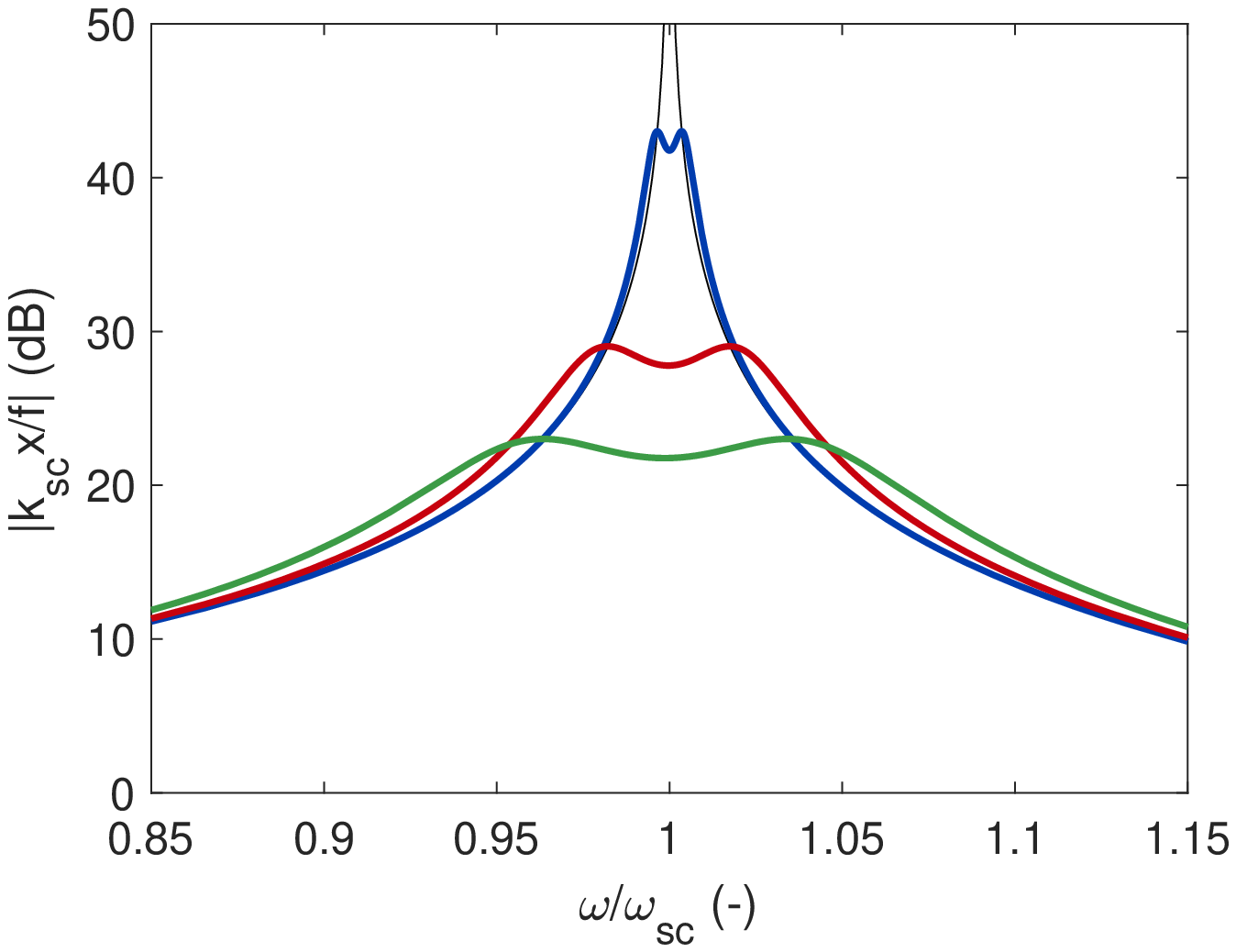}
		         \caption{}
		         \label{sfig:sdof_FRF}
		     \end{subfigure}
			\caption{Schematic representation of the electromechanical system~\subref{sfig:sdof_schematics} and FRF of the structure~\subref{sfig:sdof_FRF} when the electrodes of the transducer are short-circuited (---) and when they are connected to a parallel RL shunt circuit (\textcolor{col1}{\textbf{---}}:~$K_c=0.01$, \textcolor{col2}{\textbf{---}}:~$K_c=0.05$, \textcolor{col3}{\textbf{---}}:~$K_c=0.1$).}
			\label{fig:sdof}
		\end{figure}

		 \autoref{sfig:sdof_schematics} depicts a single-degree-of-freedom oscillator connected to a piezoelectric transducer whose electrodes are connected to a parallel RL shunt circuit, as proposed by Wu~\cite{Wu1996}. The system is governed by the following equations
		\begin{equation}
		\left\{ \begin{array}{l}
			m\ddot{x} + k_{sc}x + \gamma \dot{\psi} = f \\
			C_p^\varepsilon \ddot{\psi} + G \dot{\psi} + B \psi - \gamma \dot{x} = 0
		\end{array}\right. ,
		\label{eq:parallelRL}
		\end{equation}
		where $m$ is the structural mass, $k_{sc}=k+k_p$ is the structural stiffness with short-circuited transducer ($k$ being the stiffness of the structure without transducer, and $k_p$ the stiffness of the transducer when it is short-circuited), $\gamma$ is a coupling constant, $C_p^\varepsilon$ is the piezoelectric capacitance at constant strain, and $G=1/R$ and $B=1/L$ are the conductance (inverse of the resistance $R$) and reluctance (inverse of the inductance $L$) of the shunt circuit, respectively~\cite{Giorgio2009,Yamada2010}. $x$ is the displacement of the structure, $f$ is the external forcing, and $\psi = \int V dt$ (where $V$ is the voltage across the electrodes of the transducer) is the flux linkage. From the parameters in \autoref{eq:parallelRL}, a short-circuit resonance frequency $\omega_{sc}$ (when $\psi=0$) and a dimensionless quantity called the effective electromechanical coupling factor (EEMCF)~\cite{Thomas2009} $K_c$ can be defined as		
		\begin{equation}
			\omega_{sc}^2  =\dfrac{k_{sc}}{m}, \qquad K_c^2 = \dfrac{\gamma^2}{C_p^\varepsilon k_{sc}}.
			\label{eq:EEMCF}
		\end{equation}
		The EEMCF assesses the coupling between the mechanical and electrical dynamics. The resistance and inductance of the shunt circuit may be tuned in order to reduce the maximum amplitude of vibration of the mechanical oscillator under harmonic excitation. Generally, this requires to know the parameters $C_p^\varepsilon$, $\omega_{sc}$ and $K_c$. For instance, using the tuning rules from Yamada et al~\cite{Yamada2010},		
		\begin{equation}
			B = \dfrac{(2-K_c^2)\omega_{sc}^2 C_p^\varepsilon}{2}, \qquad G = \sqrt{\dfrac{3K_c^2}{2}}\omega_{sc}C_p^\varepsilon.
			\label{eq:Yamada}
		\end{equation}
		In general $K_c^2 \ll 1$, and thus $B/C_p^\varepsilon \approx \omega_{sc}^{2}$, i.e., the electrical resonance frequency is very close to the short-circuit mechanical one. The amplitude of the frequency response function (FRF) of the structure with this shunt circuit is shown in \autoref{sfig:sdof_FRF} for various values of the EEMCF. Clearly, the stronger the coupling, the better the performance in terms of vibration attenuation.
		
		The purpose of this paper is to propose a method able to extend this passive control strategy to potentially complex structures having multiple modes to be damped and multiple piezoelectric transducers.

\section{Dynamics of the electromechanical system}
\label{sec:dynamics}    

\begin{figure}[!ht]
        \centering
        \includegraphics[width=0.8\textwidth]{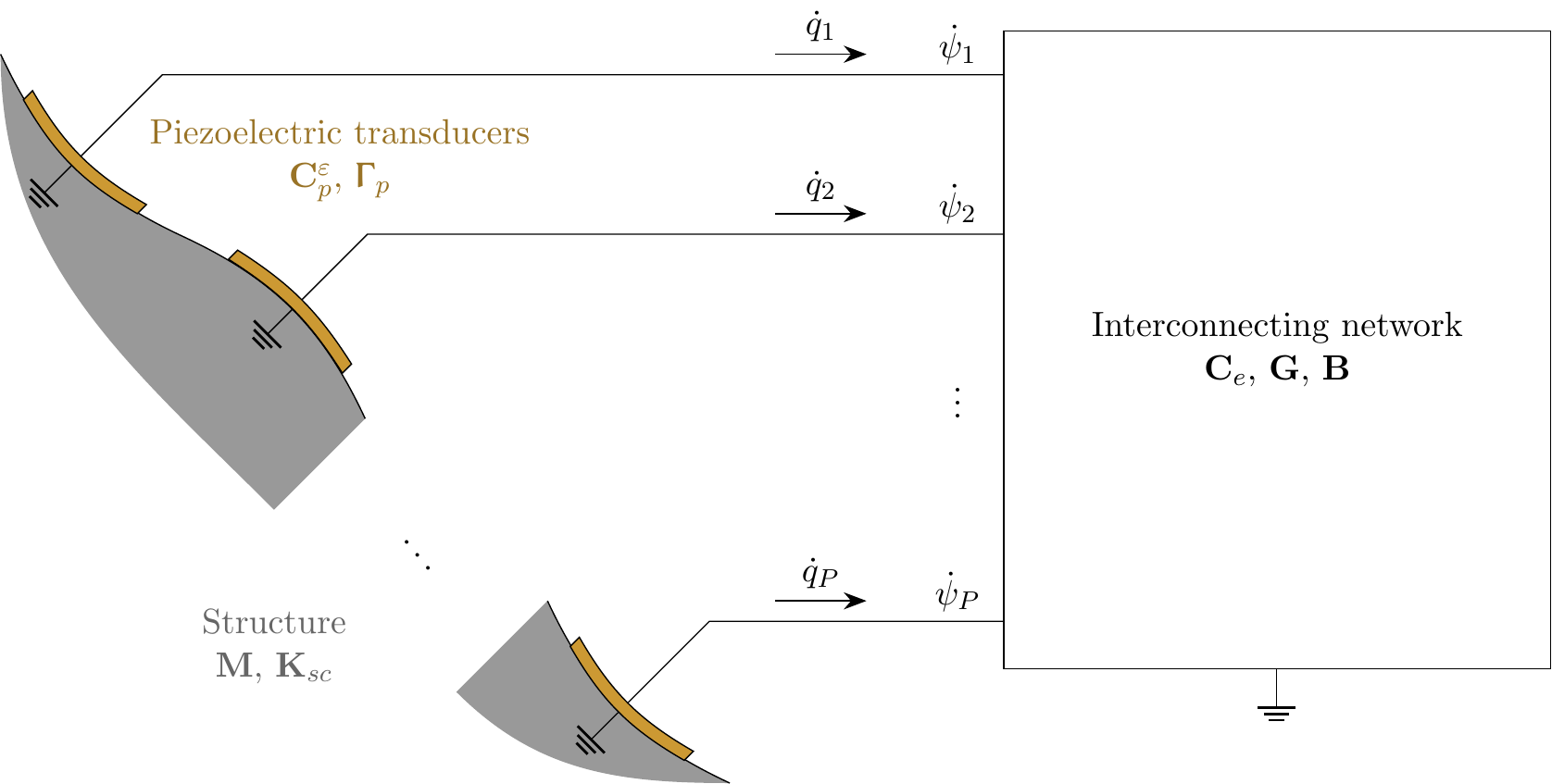}
        \caption{Schematic of a structure (in gray) with multiple piezoelectric transducers (in orange) connected to an electrical network.}
        \label{fig:generalSchematics}
    \end{figure}

    In order to achieve multimodal damping with multiple piezoelectric transducers, a properly tuned \textit{interconnecting network} can be electrically connected to these transducers, as shown in \autoref{fig:generalSchematics}. In this section, it is assumed that the network has known electrical characteristics, and an analysis of its coupling with the resonant modes of the structure is carried out. This will highlight the importance of the modal characteristics of the \textit{overall network} obtained by combining the electrical properties of the interconnecting network together with those of the transducers. These modal characteristics will then be used as design variables in the next section in order to synthesize the network.
    
    \subsection{Governing equations}
    The dynamics of a structure with $N$ degrees of freedom and with $P$ piezoelectric transducers can be described by the following system~\cite{Thomas2009}
    \begin{equation}
        \left\{
        \begin{array}{l}
            \mathbf{M} \ddot{\mathbf{x}} + \mathbf{K}_{sc} \mathbf{x} + \bm{\Upgamma}_p \dot{\bm{\uppsi}}_p = \mathbf{f} \\
            \bm{\Upgamma}_p^T \dot{\mathbf{x}} - \mathbf{C}_p^\varepsilon \ddot{\bm{\uppsi}}_p  = \dot{\mathbf{q}}_p
        \end{array}
        \right. ,
        \label{eq:piezoStructure}
    \end{equation}
    where $\mathbf{M}$ and $\mathbf{K}_{sc}$ are the structural mass and stiffness (with short-circuited transducers) matrices, $\bm{\Upgamma}_p$ is a piezoelectric coupling matrix and $\mathbf{C}_p^\varepsilon$ is the capacitance matrix of the transducers at constant strain. The vectors $\mathbf{x}$, $\mathbf{f}$, $\bm{\uppsi}_p$ and $\mathbf{q}_p$ represent the generalized mechanical degrees of freedom, the associated vector of generalized loading, the vector of flux linkages across the transducers and charges flowing out of them (considered positive if they flow out of the ungrounded electrode of the transducer), respectively.
    
    An electrical network can be electrically connected to the piezoelectric transducers in order to control multiple resonances. It is assumed that this interconnecting network is governed by the following equations
    \begin{equation}
         \mathbf{C}_e\ddot{\bm{\uppsi}} + \mathbf{G} \dot{\bm{\uppsi}} + \mathbf{B} \bm{\uppsi}  = \dot{\mathbf{q}},
        \label{eq:networkDynamics}
    \end{equation}
    where $ \mathbf{C}_e$, $\mathbf{G}$ and $\mathbf{B}$ are capacitance, conductance and reluctance matrices of the interconnecting network, respectively. The vector of flux linkages
    \begin{equation}
        \bm{\uppsi}^T = \begin{bmatrix} \bm{\uppsi}_p^T & \bm{\uppsi}_i^T \end{bmatrix},
        \label{eq:electricalDofs}
    \end{equation}
    describes the flux linkages of the $P$ ports to be connected to the piezoelectric transducers $\bm{\uppsi}_p$ and $I$ possible internal electrical degrees of freedom $\bm{\uppsi}_i$. In total, the number of these degrees of freedom is noted $N_e = P+I$. The vector $\mathbf{q}$ then describes the charges flowing into the associated ports. Connecting the network to the transducers imposes the charges flowing in the ports of the interconnecting network connected to the piezoelectric transducers to be equal to those flowing out of the transducers, as shown in \autoref{fig:generalSchematics}, leading to the relation
    \begin{equation}
        \dot{\mathbf{q}} = \mathbf{E}_p\dot{\mathbf{q}}_p
        \label{eq:networkCurrents}
    \end{equation}
    where $\mathbf{E}_p$ is a localization matrix given by
    \begin{equation}
        \mathbf{E}_p^T = \begin{bmatrix} \mathbf{I}_{P \times P} & \mathbf{0}_{P \times I} \end{bmatrix},
    \end{equation}
    $\mathbf{I}$ is the identity matrix and $\mathbf{0}$ is the zero matrix. The ports associated to internal degrees of freedom are not fed with an external current ($\dot{\mathbf{q}}_i = \mathbf{0}$). According to \autoref{eq:electricalDofs}, the flux linkages of the ports connected to the piezoelectric transducers are also given by
    \begin{equation}
        \bm{\uppsi}_p = \mathbf{E}_p^T\bm{\uppsi}.
        \label{eq:networkFluxLinkage}
    \end{equation}
    Upon connecting the network governed by \autoref{eq:networkDynamics} to the transducers of the piezoelectric structure, the governing equations in \autoref{eq:piezoStructure} become, using Equations~\ref{eq:networkCurrents} and \ref{eq:networkFluxLinkage},
    \begin{equation}
        \left\{
        \begin{array}{l}
            \mathbf{M} \ddot{\mathbf{x}} + \mathbf{K}_{sc} \mathbf{x} + \bm{\Upgamma}_p\mathbf{E}_p^T \dot{\bm{\uppsi}} = \mathbf{f} \\
            \left( \mathbf{C}_e + \mathbf{E}_p\mathbf{C}_p^\varepsilon\mathbf{E}_p^T \right)\ddot{\bm{\uppsi}} + \mathbf{G} \dot{\bm{\uppsi}} + \mathbf{B} \bm{\uppsi} - \mathbf{E}_p\bm{\Upgamma}_p^T \dot{\mathbf{x}} = \mathbf{0}
        \end{array}
        \right. ,
        \label{eq:electromechanicalSystem}
    \end{equation}
    which describe the dynamics of the structure coupled to the electrical network. Specifically, the second set of equations describes the dynamics of the overall network. In the sequel, the capacitance matrix of the overall network will be noted $\mathbf{C}$ for brevity, i.e. 
    \begin{equation}
        \mathbf{C} =  \mathbf{C}_e + \mathbf{E}_p\mathbf{C}_p^\varepsilon\mathbf{E}_p^T.
        \label{eq:overallCapacitanceMatrix}
    \end{equation}

    \subsection{Mechanical and electrical modes}
        
        The dynamics of the uncoupled ($\bm{\Upgamma}_p=\mathbf{0}$) mechanical and electrical systems are considered herein.
        
        The mechanical modes with short-circuited transducers are defined by the generalized eigenvalue problem~\cite{Geradin2014}
        \begin{equation}
            \mathbf{K}_{sc}\bm{\Upphi}_{sc} = \mathbf{M}\bm{\Upphi}_{sc}\bm{\Upomega}^2_{sc}, \qquad \bm{\Upomega}^2_{sc} = \begin{bmatrix} \omega_{sc,1}^2 & & \\ & \ddots & \\ & & \omega_{sc,N}^2\end{bmatrix},
        \end{equation}
        where $\bm{\Upphi}_{sc}$ is the mechanical mode shapes matrix and $\bm{\Upomega}^2_{sc}$ is a diagonal matrix whose entries are the squared natural frequencies of the structure. The modes verify a series of orthogonality relationships, and are usually mass-normalized, i.e.,
        \begin{equation}
            \bm{\Upphi}^T_{sc} \mathbf{M} \bm{\Upphi}_{sc} = \mathbf{I}, \qquad \bm{\Upphi}_{sc}^T\mathbf{K}_{sc}\bm{\Upphi}_{sc} = \bm{\Upomega}^2_{sc}.
            \label{eq:mechanicalModesOrthogonality}
        \end{equation}
        
        In a similar way, electrical modes of the lossless network ($\mathbf{G} = \mathbf{0}$) can be defined by the following equation:
        \begin{equation}
            \mathbf{B}\bm{\Upphi}_e = \mathbf{C}\bm{\Upphi}_e\bm{\Upomega}^2_e, \qquad \bm{\Upomega}^2_{e} = \begin{bmatrix} \omega_{e,1}^2 & & \\ & \ddots & \\ & & \omega_{e,N_e}^2\end{bmatrix},
        \end{equation}
        where $\bm{\Upphi}_e$ is the electrical mode shapes matrix and $\bm{\Upomega}_e^2$ is a diagonal matrix whose entries are the squared natural frequencies of the network. The mode shapes also verify orthogonality properties and can be capacitance-normalized
        \begin{equation}
            \bm{\Upphi}^T_e \mathbf{C} \bm{\Upphi}_e = \mathbf{I}, \qquad \bm{\Upphi}^T_e\mathbf{B}\bm{\Upphi}_e = \bm{\Upomega}^2_e.
            \label{eq:electricalModesOrthogonality}
        \end{equation}
        It is also assumed that the electrical conductance matrix is such that the mode shapes also respect the same kind of orthogonality relations
        \begin{equation}
            \bm{\Upphi}^T_e\mathbf{G}\bm{\Upphi}_e = 2\mathbf{Z}_e\bm{\Upomega}_e, \qquad \mathbf{Z}_e = \begin{bmatrix} \zeta_{e,1} & & \\ & \ddots & \\ & & \zeta_{e,N_e}\end{bmatrix}
            \label{eq:electricalModesOrthogonality2},
        \end{equation}
        where $\mathbf{Z}_e$ is a diagonal matrix containing the electrical modal damping ratios.
    
    \subsection{Modal coupling}
        Assuming that the $k^{th}$ electrical resonance frequency is close to the $r^{th}$ mechanical one, and for forcing frequencies near them, the dynamics of the mechanical and electrical systems are assumed to be dominated by their associated uncoupled modes, i.e., 
        \begin{equation}
            \mathbf{x} \approx \bm{\upphi}_{sc,r} \eta_{sc,r}, \qquad \bm{\uppsi} \approx \bm{\upphi}_{e,k}\eta_{e,k},
            \label{eq:oneModeHypothesis}
        \end{equation}
        where $\eta_{sc,r}$ and $\eta_{e,k}$ are mechanical and electrical modal coordinates, respectively, $\bm{\upphi}_{sc,r}$ is the $r^{th}$ column of $\bm{\Upphi}_{sc}$ and $\bm{\upphi}_{e,k}$ is the $k^{th}$ column of $\bm{\Upphi}_e$. Inserting this ansatz into \autoref{eq:electromechanicalSystem}, premultiplying them by $\bm{\upphi}_{sc,r}^T$ and $\bm{\upphi}_{e,k}^T$ and accounting for Equations~\ref{eq:mechanicalModesOrthogonality}, \ref{eq:electricalModesOrthogonality} and \ref{eq:electricalModesOrthogonality2}, one gets the reduced system
        \begin{equation}
            \left\{
            \begin{array}{l}
                 \ddot{\eta}_{sc,r} + \omega_{sc,r}^2\eta_{sc,r} + \gamma_{rk} \dot{\eta}_{e,k} = \bm{\upphi}_{sc,r}^T\mathbf{f} \\
                 \ddot{\eta}_{e,k} + 2\zeta_{e,k}\omega_{e,k}\dot{\eta}_{e,k} + \omega_{e,k}^2\eta_{e,k} - \gamma_{rk}\dot{\eta}_{sc,r} = 0
            \end{array}
            \right. ,
            \label{eq:oneModeSystem}
        \end{equation}
        where 
        \begin{equation}
            \gamma_{rk} = \bm{\upphi}_{sc,r}^T\bm{\Upgamma}_p\mathbf{E}_p^T\bm{\upphi}_{e,k}
            \label{eq:modalECF}
        \end{equation}
        is a modal electromechanical coupling coefficient between mechanical mode $r$ and electrical mode $k$. \autoref{eq:oneModeHypothesis} is a central assumption of this paper. It can be used to accurately describe the uncoupled systems in a certain frequency range, and the coupled system can be well-approximated by retaining in \autoref{eq:oneModeSystem} the modal coupling coefficient defined in \autoref{eq:modalECF}. A similar hypothesis was made in~\cite{Alessandroni2002} using the eigenfunctions of the spatial differential operators of the uncoupled systems. Noting the similarity between \autoref{eq:parallelRL} and \autoref{eq:oneModeSystem}, an EEMCF can be defined by analogy to \autoref{eq:EEMCF} as
        \begin{equation}
            \widehat{K}_{c,rk}^2 = \dfrac{\gamma_{rk}^2}{\omega_{sc,r}^2}.
            \label{eq:modalEEMCF}
        \end{equation}

\section{Design of an electrical network}

    The goal of this section is to find optimal electrical modal parameters for the network and to deduce its electrical capacitance, conductance and reluctance matrices from Equations~\ref{eq:electricalModesOrthogonality} and \ref{eq:electricalModesOrthogonality2}. Similarly to what has been explained in \autoref{sec:piezoShunt}, the natural frequencies of the network should be close to those of a targeted set of structural modes~\cite{Vidoli2000}. The choice of the electrical mode shapes is less straightforward. As in \autoref{sec:piezoShunt}, it is desirable from a vibration reduction performance perspective to maximize the EEMCF given in \autoref{eq:modalEEMCF}. With \autoref{eq:modalECF}, it is thus sought to maximize the amplitude of the electrical mode shapes at the piezoelectric transducers. Their amplitude is however limited by passivity constraints, as shall be seen.

    \label{sec:networkDesign}
    \subsection{Passivity}

        The passivity of the network is set as a design requirement in this work. Passive control has the advantage of guaranteeing the stability of the controlled system (because a passive system can only store or dissipate energy), and theoretically does not need external power sources for operation. It will be shown that this requirement places limits on the attainable modal characteristics, which in turn limits the performance of the passive control system. The electrical modal characteristics can be optimized under the passivity constraints, and these optimal characteristics can finally be used to synthesize the electrical network.
    
        Using the Laplace transform of \autoref{eq:electromechanicalSystem}, the governing equations can also be expressed in the frequency domain by 
        \begin{equation}
            \left\{
            \begin{array}{l}
                \mathbf{H}(s)\mathbf{x} + s\bm{\Upgamma}_p\mathbf{E}_p^T \bm{\uppsi} = \mathbf{f} \\
                \mathbf{Y}(s) \bm{\uppsi} - \mathbf{E}_p\bm{\Upgamma}_p^T \mathbf{x} = \mathbf{0}
            \end{array}
            \right. ,
            \label{eq:electromechanicalSystemLaplace}
        \end{equation}
        where $s$ is the Laplace variable,
        \begin{equation}
            \mathbf{H}(s) = \mathbf{M}s^2 + \mathbf{K}_{sc}
        \end{equation}
        is the dynamic stiffness matrix~\cite{Geradin2014}, and
        \begin{equation}
            \mathbf{Y}(s) = s \mathbf{C} + \mathbf{G} + \dfrac{1}{s}\mathbf{B}
            \label{eq:admittanceMatrix}
        \end{equation}
        is the nodal admittance matrix~\cite{Grainger1994}.
        
        According to Gannett and Chua~\cite{Gannett1978}, the nodal admittance matrix given in \autoref{eq:admittanceMatrix} must fulfill the conditions
		\begin{enumerate}[label=(\roman*)]
			\item $\mathbf{Y}(s)$ has no poles in $\left\{s\in \mathbb{C} | \Re(s) > 0\right\}$ ($\Re$ denotes the real part operator), \label{cond:i}
			\item $\mathbf{Y}(\sigma)$ is a real matrix for $\sigma \in \mathbb{R}^+ $,\label{cond:ii}
			\item $\mathbf{Y}(s) + \mathbf{Y}^H(s)$ is positive semidefinite in $\left\{s\in \mathbb{C} | \Re(s) > 0\right\}$ (superscript $H$ denotes Hermitian transposition),\label{cond:iii}
			\item The network associated to $\mathbf{Y}$ is controllable\label{cond:iv}
		\end{enumerate}		
		in order to be the admittance matrix of a passive network (i.e., realizable using passive capacitors, resistors, inductors and ideal transformers). According to \autoref{eq:admittanceMatrix}, $\mathbf{Y}(s)$ has one simple pole at $s=0$ (for a nonzero reluctance matrix), so Condition~\ref{cond:i} is satisfied. Condition~\ref{cond:ii} is satisfied since $\mathbf{C}$, $\mathbf{G}$ and $\mathbf{B}$ are real matrices. Condition~\ref{cond:iv} is also verified, because all the electrical states are controllable and observable for an admittance matrix of the form of \autoref{eq:admittanceMatrix}. Now, since $\mathbf{C}$, $\mathbf{G}$ and $\mathbf{B}$ are also symmetric matrices, using \autoref{eq:admittanceMatrix}, Condition~\ref{cond:iii} becomes
		\begin{equation}
			\mathbf{Y}(\sigma + j\omega) + \mathbf{Y}^H(\sigma + j\omega)  = 2 \sigma \mathbf{C} + 2\mathbf{G} + \dfrac{2\sigma}{\sigma^2 + \omega^2}\mathbf{B} \succeq 0,
		\end{equation}
		where $\succeq 0$ indicates that the matrix is positive semidefinite. For $\sigma > 0$ and $\omega\in \mathbb{R}$, this matrix is positive semidefinite if $\mathbf{C}$, $\mathbf{G}$ and $\mathbf{B}$ are positive semidefinite themselves, which gives the criteria to satisfy Condition~\ref{cond:iii}. Because of Equations~\ref{eq:electricalModesOrthogonality} and \ref{eq:electricalModesOrthogonality2}, the matrices $\mathbf{C}$, $\mathbf{G}$ and $\mathbf{B}$ are guaranteed to be positive semidefinite, which eventually ensures the passivity of the overall network.
		
		Now, one must consider that the piezoelectric transducers are integrated into the overall network associated with the matrix $\mathbf{C}$. The capacitance matrix $ \mathbf{C}_e$ of the interconnecting network is obtained by removing the contribution of the piezoelectric capacitance from $\mathbf{C}$, as indicated by \autoref{eq:overallCapacitanceMatrix}. Thus, if $ \mathbf{C}_e$ is the capacitance matrix of a passive interconnecting network, it must satisfy
		\begin{equation}
		     \mathbf{C}_e = \mathbf{C} - \mathbf{E}_p\mathbf{C}_p^\varepsilon\mathbf{E}_p^T \succeq 0.
		    \label{eq:posDefCn}
		\end{equation}
        A positive semidefinite matrix has positive eigenvalues, and a necessary (but not sufficient) condition for this is that the determinant of this matrix (being equal to the product of the eigenvalues) must be positive. This determinant is given by
        \begin{equation}
            \begin{array}{rl}
		    \det \left( \mathbf{C} - \mathbf{E}_p\mathbf{C}_p^\varepsilon\mathbf{E}_p^T \right) &= \det\left(\mathbf{C}^{1/2}\right)\det \left( \mathbf{I} - \mathbf{C}^{-1/2}\mathbf{E}_p\mathbf{C}_p^\varepsilon\mathbf{E}_p^T\mathbf{C}^{-1/2} \right)\det\left(\mathbf{C}^{1/2}\right)  \\ &= \det \left( \mathbf{I} - \left(\mathbf{C}_p^\varepsilon\right)^{1/2}\mathbf{E}_p^T\mathbf{C}^{-1}\mathbf{E}_p\left(\mathbf{C}_p^\varepsilon\right)^{1/2} \right)\det\left(\mathbf{C}\right) \\
		    &= \det \left( \left(\mathbf{C}_p^\varepsilon\right)^{-1} - \mathbf{E}_p^T\mathbf{C}^{-1}\mathbf{E}_p \right)\det\left(\mathbf{C}\right)\det\left(\mathbf{C}_p^\varepsilon\right)
		    \end{array},
		\end{equation}
        where the fact that the determinant of a product of matrices is the product of their determinants has been used, as well as the Weinstein-Aronszajn formula (\cite{Chervov2009}, Equation (5.47)). Since both $\mathbf{C}$ and $\mathbf{C}_p^\varepsilon$ are positive definite, a necessary condition for the positive semidefiniteness of $ \mathbf{C}_e$ is thus
        \begin{equation}
            \left( \mathbf{C}_p^\varepsilon \right)^{-1} - \mathbf{E}_p^T\bm{\Upphi}_e\bm{\Upphi}_e^T\mathbf{E}_p \succeq 0.
            \label{eq:passivityCondition}
        \end{equation}
        If the modal amplitudes $\bm{\Upphi}_e^T\mathbf{E}_p$ are initially zero, this condition will be satisfied because $\mathbf{C}_p^\varepsilon$ is itself positive definite. Increasing gradually these modal amplitudes will affect continuously the eigenvalues of the matrix in \autoref{eq:passivityCondition}, up to the point where one of them reaches zero. Increasing the amplitudes beyond this point, the matrix will no longer be positive semidefinite, and the interconnecting network will no longer be passive.

    \subsection{Optimal electrical mode shapes}
        
        \autoref{eq:passivityCondition} constrains the amplitude of the electrical mode shapes but leaves freedom on their actual shape. However, some choices yield better electromechanical coupling factors than others. The purpose of this subsection is to provide an expression for the set of electrical mode shapes maximizing the modal electromechanical coupling coefficients.
        
        \subsubsection{Optimization for a single mode}
        
        It is first considered that the network targets a single mode of the structure. To simplify the exposition, the electrical mode shape at the piezoelectric transducers is made dimensionless by the transformation 
        \begin{equation}
            \bm{\upvarphi}_p = \left(\mathbf{C}_p^\varepsilon\right)^{1/2} \mathbf{E}_p^T \bm{\upphi}_e = \left(\mathbf{C}_p^\varepsilon\right)^{1/2} \bm{\upphi}_p.
            \label{eq:dimensionlessModeShape}
        \end{equation}
        Inserting this relation into \autoref{eq:passivityCondition} yields
        \begin{equation}
            \left( \mathbf{C}_p^\varepsilon \right)^{-1} - \left(\mathbf{C}_p^\varepsilon\right)^{-1/2}\bm{\upvarphi}_p\bm{\upvarphi}_p^T\left(\mathbf{C}_p^\varepsilon\right)^{-1/2} \succeq 0.
        \end{equation}
        The positive semidefinite character is not altered by a pre- and postmultiplication by $\left(\mathbf{C}_p^\varepsilon\right)^{1/2}$, thus, 
        \begin{equation}
            \mathbf{I} - \bm{\upvarphi}_{p}\bm{\upvarphi}_{p}^T \succeq 0.
            \label{eq:passivityConstraintOneMode}
        \end{equation}
        It is useful to point out here that if $\mathbf{A}$ is any $n\times n$ matrix with eigenvalues $\lambda_1, \cdots , \lambda_n$, the eigenvalues of $\mathbf{I}-\mathbf{A}$ are $1-\lambda_1, \cdots , 1-\lambda_n$. \autoref{eq:passivityConstraintOneMode} is therefore satisfied if the dimensionless mode shape satisfies
        \begin{equation}
            \left| \bm{\upvarphi}_{p} \right| \leq 1.
        \end{equation}
        
        Now, the EEMCF defined in \autoref{eq:modalEEMCF} is maximized if the modal electromechanical coupling coefficient defined in \autoref{eq:modalECF} is itself maximized, leading to the following constrained optimization problem
        \begin{equation}
            \begin{array}{rl}
                \underset{\bm{\upvarphi}_p}{\text{Maximize}} & \bm{\upphi}_{sc,r}^T\bm{\Upgamma}_p\left(\mathbf{C}_p^\varepsilon\right)^{-1/2}\bm{\upvarphi}_{p} \\
                \text{Subject to} & \bm{\upvarphi}_{p}^T\bm{\upvarphi}_{p} \leq 1
            \end{array}.
        \end{equation}
        The admissible solution to this optimization problem can be found as
        \begin{equation}
            \bm{\upvarphi}_{p}^\star = \dfrac{1}{\sqrt{\bm{\upphi}_{sc,r}^T\bm{\Upgamma}_p\left(\mathbf{C}_p^\varepsilon\right)^{-1}\bm{\Upgamma}_p^T\bm{\upphi}_{sc,r}}}\left(\mathbf{C}_p^\varepsilon\right)^{-1/2}\bm{\Upgamma}_p^T\bm{\upphi}_{sc,r},
        \end{equation}
        and the associated optimal capacitance-normalized mode shape at the piezoelectric transducers is retrieved with \autoref{eq:dimensionlessModeShape} as
        \begin{equation}
            \bm{\upphi}_{p}^\star =  \dfrac{1}{\sqrt{\bm{\upphi}_{sc,r}^T\bm{\Upgamma}_p\left(\mathbf{C}_p^\varepsilon\right)^{-1}\bm{\Upgamma}_p^T\bm{\upphi}_{sc,r}}}\left(\mathbf{C}_p^\varepsilon\right)^{-1}\bm{\Upgamma}_p^T\bm{\upphi}_{sc,r}.
            \label{eq:optimalElectricalModeShape}
        \end{equation}
        
        \subsubsection{Optimization for multiple modes}
            In order for the network to mitigate $N_s$ structural resonances, it is assumed that the dimensionless electrical mode shapes are simply scaled versions of their optimal counterparts in the single-mode case, i.e.,
            \begin{equation}
                \overline{\bm{\Upphi}}_p = \begin{bmatrix} \bm{\upvarphi}_{p,1} & \cdots & \bm{\upvarphi}_{p,N_s}  \end{bmatrix} = \overline{\bm{\Upphi}}_p^\star \mathbf{D}_p = \begin{bmatrix} \bm{\upvarphi}_{p,1}^\star & \cdots & \bm{\upvarphi}_{p,N_s}^\star \end{bmatrix}\begin{bmatrix}d_{p,1} & & \\ & \ddots & \\ & & d_{p,N_s} \end{bmatrix},
            \end{equation}
            where $d_{p,k}$ is the positive scaling factor associated to electrical mode $k$. In order to satisfy the passivity constraints (\autoref{eq:passivityCondition}), these scaling factors must be chosen such that
            \begin{equation}
                \mathbf{I} - \overline{\bm{\Upphi}}_p^\star \mathbf{D}_p^2 \left(\overline{\bm{\Upphi}}_p^\star\right)^T \succeq 0.
                \label{eq:passivityConstraintModal}
            \end{equation}
            To determine them, a set of positive relative scaling factors $\overline{d}_{p,k}$ can be chosen arbitrarily. Their magnitude does not matter, but their relative magnitude can be set so as to put more control authority on specific modes at the expense of control authority on other modes. The actual scaling factors are obtained by multiplication by a scalar factor $\alpha$
            \begin{equation}
                \mathbf{D}_p = \alpha \begin{bmatrix}\overline{d}_{p,1} & & \\ & \ddots & \\ & & \overline{d}_{p,N_s} \end{bmatrix} = \alpha \overline{\mathbf{D}}_p,
            \end{equation}
            which can be determined to enforce the passivity constraints. From \autoref{eq:passivityConstraintModal}, it should satisfy
            \begin{equation}
                \alpha \leq \dfrac{1}{\sqrt{\lambda_{Max} \left(\overline{\bm{\Upphi}}_p^\star \overline{\mathbf{D}}_p^2 \left(\overline{\bm{\Upphi}}_p^\star\right)^T\right)} },
                \label{eq:alphaScaling}
            \end{equation}
            where $\lambda_{Max}(\cdot)$ is the maximum eigenvalue of the matrix at hand. The capacitance-normalized mode shapes at the piezoelectric transducers can finally be retrieved as
            \begin{equation}
                \bm{\Upphi}_p = \left(\mathbf{C}_p^\varepsilon\right)^{-1/2}\overline{\bm{\Upphi}}_p = \alpha \left(\mathbf{C}_p^\varepsilon\right)^{-1/2}\overline{\bm{\Upphi}}_p^\star \overline{\mathbf{D}}_p.
                \label{eq:optimalMultipleModeShapes}
            \end{equation}
            The mode shapes obtained this way are proportional to $\alpha$. Equations~\eqref{eq:modalECF} and~\eqref{eq:modalEEMCF} then imply that the EEMCF itself is proportional to $\alpha$. Hence, the inequality in \autoref{eq:alphaScaling} should ideally be an equality to maximize the amplitude of the electrical modes at the piezoelectric transducers and thus the coupling.

            It should be pointed out that if the piezoelectric patches are all identical, $\mathbf{C}_p^\varepsilon$ is proportional to the identity matrix.  $\bm{\Upgamma}_p^T\bm{\Upphi}_{sc}$ can be seen as modal strains in the transducers, and the electrical mode shapes $\bm{\Upphi}_p$ are proportional to them, according to \autoref{eq:optimalElectricalModeShape} in this case. Hence, the network is an analog of the structure, because it features identical resonance frequencies and mode shapes, as in~\cite{Alessandroni2002} (but the scaling factors make the electrical mode shapes in \autoref{eq:optimalMultipleModeShapes} a scaled version of their mechanical counterparts).
            
    \subsection{Optimal electrical frequencies and damping ratios}
    
        In order to fully specify the characteristics of the network, its natural frequencies and associated damping ratios should also be prescribed. These modal characteristics can simply be tuned by using \autoref{eq:Yamada} with the characteristics of the system in \autoref{eq:oneModeSystem}, i.e., by replacing $C_p^\varepsilon$ by one, $\omega_{sc}$ by $\omega_{sc,r}$, $K_c$ by $\widehat{K}_{c,rk}$, $B$ by $\omega_{e,k}^2$ and $G$ by $2\zeta_{e,k}\omega_{e,k}$. The optimal frequency and damping ratio for electrical mode $k$ are thus given by
		\begin{equation}
			\omega_{e,k}^2 = \dfrac{(2-\widehat{K}_{c,rk}^2)\omega_{sc,r}^2}{2}, \qquad \zeta_{e,k} = \dfrac{\sqrt{3}}{2}\sqrt{\dfrac{\widehat{K}_{c,rk}^2}{2-\widehat{K}_{c,rk}^2}},
			\label{eq:YamadaModal}
		\end{equation}
        respectively.
        
    \subsection{Network synthesis}
    
        Now that all the modal characteristics are specified, the capacitance, conductance and reluctance matrices of the network can be computed so as to satisfy Equations~\ref{eq:electricalModesOrthogonality} and \ref{eq:electricalModesOrthogonality2}. Depending on the number of transducers $P$ compared to the number of targeted modes $N_s$, different cases must be considered.
    
        \subsubsection{$P = N_s$}
        
            The case $P=N_s$ is the simplest, because the matrix $\bm{\Upphi}_e=\bm{\Upphi}_p$ is square and presumably non-singular. The electrical matrices can then be obtained from the inversion of Equations~\ref{eq:electricalModesOrthogonality} and \ref{eq:electricalModesOrthogonality2} 
            \begin{equation}
                \mathbf{C} = \bm{\Upphi}_e^{-T}\bm{\Upphi}_e^{-1}, \qquad \mathbf{G} = 2\bm{\Upphi}_e^{-T}\mathbf{Z}_e\bm{\Upomega}_e\bm{\Upphi}_e^{-1}, \qquad \mathbf{B} = \bm{\Upphi}_e^{-T} \bm{\Upomega}^2_e \bm{\Upphi}_e^{-1}.
                \label{eq:electricalSpectralExpansion}
            \end{equation}
        
        \subsubsection{$P < N_s$}
            
            If $P<N_s$, there are less transducers than targeted structural modes, and the matrix $\bm{\Upphi}_p$ has more columns than rows. In that case, internal degrees of freedom $\bm{\uppsi}_i$ should be introduced, which also implies that the mode shapes on these degrees of freedom $\bm{\Upphi}_i$ should be specified in order to obtain a square mode shape matrix
            \begin{equation}
                \bm{\Upphi}_e = \begin{bmatrix} \bm{\Upphi}_p \\ \bm{\Upphi}_i \end{bmatrix},
            \end{equation}
            and the electrical matrices can then be obtained with \autoref{eq:electricalSpectralExpansion}. As long as the electrical mode shape matrix is non-singular, the mode shapes on the internal degrees of freedom can be chosen arbitrarily without affecting the performance of the network. 
            
        \subsubsection{$P > N_s$}
        
            If $P>N_s$, there are more transducers than targeted structural modes, and the matrix $\bm{\Upphi}_e=\bm{\Upphi}_p$ has more rows than columns. Compliance with Equations~\ref{eq:electricalModesOrthogonality} and \ref{eq:electricalModesOrthogonality2} leads to an underconstrained problem. This leaves some freedom on the choice of the other electrical mode shapes of the network (while keeping in mind the passivity requirements). Another simple way to enforce Equations~\ref{eq:electricalModesOrthogonality} and \ref{eq:electricalModesOrthogonality2} is to choose
            \begin{multline}
                \mathbf{C} = \bm{\Upphi}_e \left(\bm{\Upphi}_e^T\bm{\Upphi}_e\right)^{-2}\bm{\Upphi}_e^T + \mathbf{V}\mathbf{D}_V\mathbf{V}^T, \qquad \mathbf{G} = \bm{\Upphi}_e \left(\bm{\Upphi}_e^T\bm{\Upphi}_e\right)^{-1}2\mathbf{Z}_e\bm{\Upomega}_e\left(\bm{\Upphi}_e^T\bm{\Upphi}_e\right)^{-1}\bm{\Upphi}_e^T, \\ \mathbf{B} = \bm{\Upphi}_e \left(\bm{\Upphi}_e^T\bm{\Upphi}_e\right)^{-1}\bm{\Upomega}_e^2\left(\bm{\Upphi}_e^T\bm{\Upphi}_e\right)^{-1}\bm{\Upphi}_e^T, 
                \label{eq:electricalMatricesCase3}
            \end{multline}
            where $\mathbf{D}_V$ is a diagonal matrix, and $\mathbf{V}$ contains the orthogonal basis of the kernel of $\bm{\Upphi}_e^T$, i.e.,
            \begin{equation}
                \bm{\Upphi}_e^T\mathbf{V} = \mathbf{0}, \qquad \mathbf{V}^T\mathbf{V} = \mathbf{I}.
            \end{equation}
            $\mathbf{D}_V$ does not affect the electromechanical coupling with the targeted modes, but has to be chosen such that \autoref{eq:posDefCn} is satisfied. The projection of this equation on the subspace spanned by $\bm{\Upphi}_e$ satisfies that equation, and that on $\mathbf{V}$ yields the condition
            \begin{equation}
                \mathbf{D}_V - \mathbf{V}^T\mathbf{E}_p\mathbf{C}_p^\varepsilon\mathbf{E}_p^T\mathbf{V} \succeq 0
            \end{equation}
            choosing 
            \begin{equation}
                \mathbf{D}_V = \beta \mathbf{I}
            \end{equation}
            for simplicity, \autoref{eq:posDefCn} can be satisfied if 
            \begin{equation}
                \beta = \lambda_{Max}\left(\mathbf{V}^T\mathbf{E}_p\mathbf{C}_p^\varepsilon\mathbf{E}_p^T\mathbf{V}\right).
            \end{equation}
    \subsection{Design procedure}
        
\begin{figure}[!ht]
            \centering
            \includegraphics[width=\textwidth]{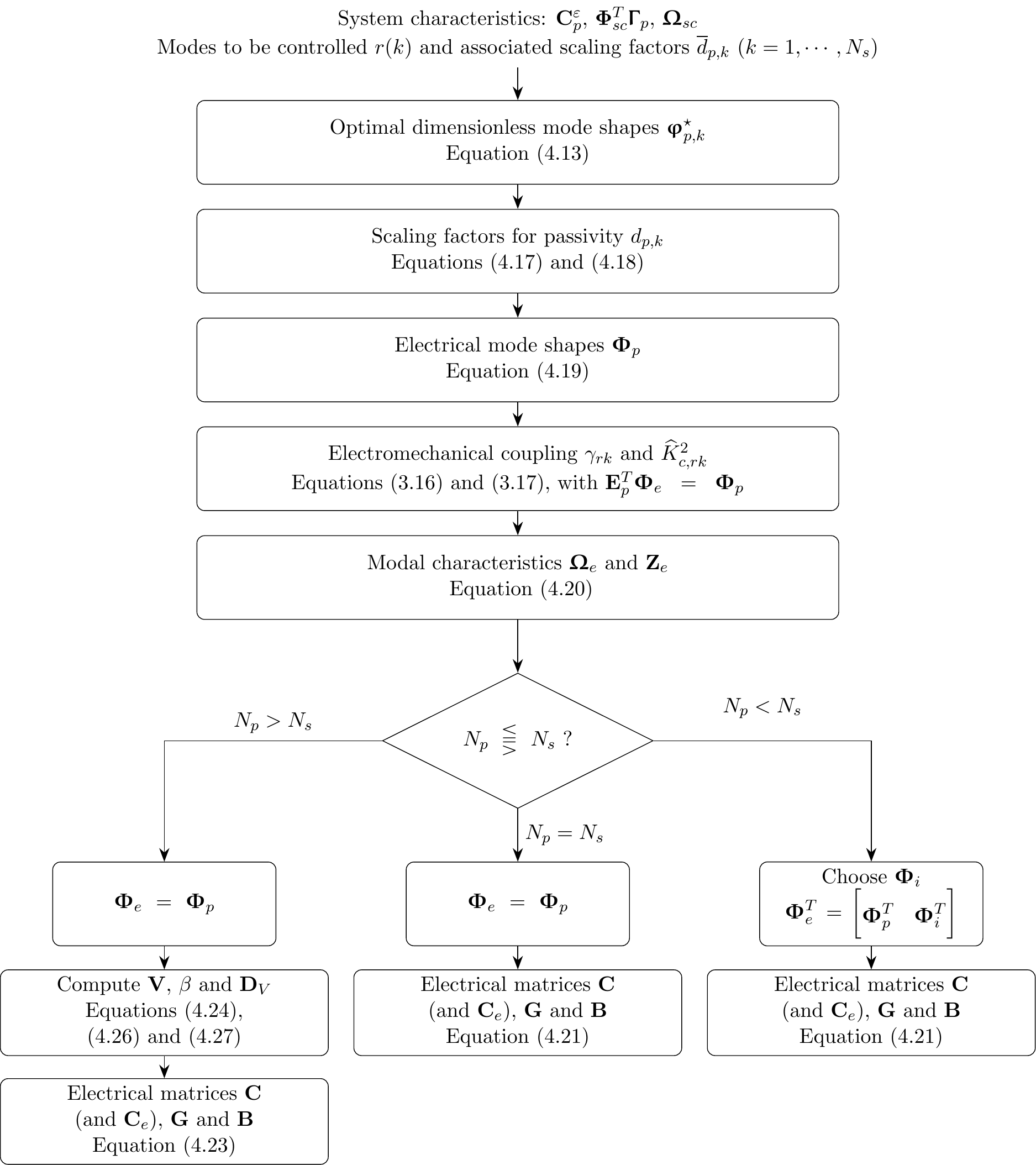}
            \caption{Flowchart of the proposed modal-based synthesis.}
            \label{fig:flowchart}
        \end{figure}

        The proposed design procedure is summarized in \autoref{fig:flowchart}: from a set of targeted modes and associated scaling factors reflecting the desired control authority on these modes, the electrical matrices of the network can be obtained. Given the central role played by modal properties in this approach, it is named modal-based synthesis.

\section{Numerical examples}
\label{sec:examples}

    \subsection{Free-free beam}

        The free-free beam in~\cite{Lossouarn2015a} is the first example used to demonstrate the modal-based approach and to compare it to methods based on electrical analogs. The beam is composed of 20 identical cells where a pair of thin piezoelectric patches symmetrically bonded to the structure and electrically connected in parallel can be used to damp the vibrations of the beam. The beam was modeled using the finite element method~\cite{Thomas2009} with ten elements per cell. This yielded the matrices $\mathbf{M}$, $\mathbf{K}_{sc}$ and $\bm{\Upgamma}_p$.
        
        \subsubsection{Damping of the first four modes}
        
\begin{figure}[!ht]
            \centering
            \includegraphics[width=0.8\textwidth]{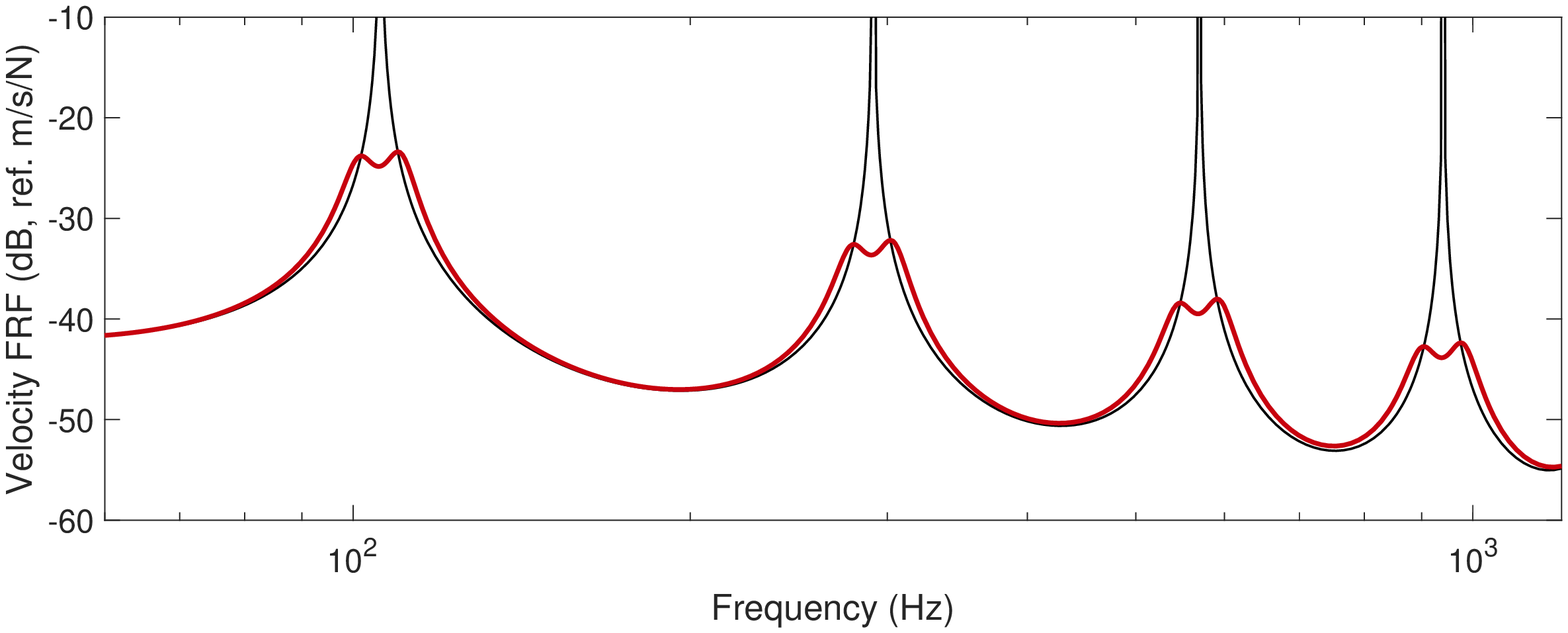}
            \caption{Velocity FRF of the beam with short-circuited patches~(---), with the network synthesized with the modal-based approach~(\textbf{\textcolor{col2}{---}}).}
            \label{fig:beamFRF0}
        \end{figure}

        A network targeting the first four flexible modes of the beam was synthesized with the modal-based method described in Figure \ref{fig:flowchart}, choosing unit scaling factors for each mode ($\overline{d}_{p,k}=1$, $k=1,2,3,4$). The FRF of the beam transversely excited at one end and whose transverse displacement is measured at the other end is shown in \autoref{fig:beamFRF0}. It is observed that the electrical network can very effectively damp out the resonant vibrations of the four targeted modes.

        \subsubsection{Relative scaling factors}
        
\begin{figure}[!ht]
            \centering
            \includegraphics[width=0.8\textwidth]{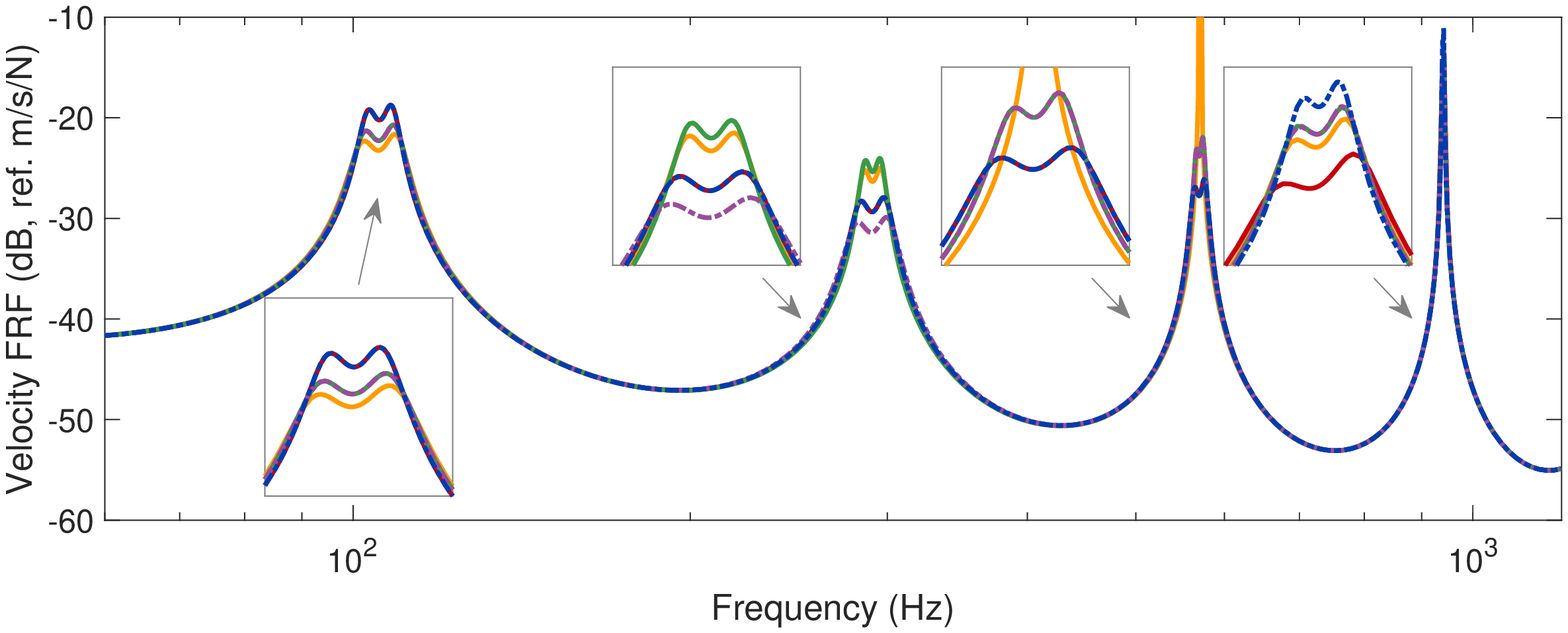}
            \caption{Velocity FRF of the beam with grouped patches connected to a network synthesized with the modal-based approach: $[\overline{d}_{p,1}, \overline{d}_{p,2} ,\overline{d}_{p,3} ,\overline{d}_{p,4}]= [1,1,1,1]$ (\textbf{\textcolor{col2}{---}}), $[\overline{d}_{p,1}, \overline{d}_{p,2} ,\overline{d}_{p,3} ,\overline{d}_{p,4}]= [2,1,1,1]$ (\textbf{\textcolor{col3}{---}}), $[\overline{d}_{p,1}, \overline{d}_{p,2} ,\overline{d}_{p,3} ,\overline{d}_{p,4}]= [2,1,0,1]$ (\textbf{\textcolor{col5}{---}}), $[\overline{d}_{p,1}, \overline{d}_{p,2} ,\overline{d}_{p,3} ,\overline{d}_{p,4}]= [2,2,1,1]$ (\textbf{\textcolor{col4}{-$\cdot$-}}) and $[\overline{d}_{p,1}, \overline{d}_{p,2} ,\overline{d}_{p,3} ,\overline{d}_{p,4}]= [2,2,2,1]$ (\textbf{\textcolor{col1}{-$\cdot$-}}).}
            \label{fig:beamFRFScaling}
        \end{figure}

\begin{figure}
            \centering
            \includegraphics[width=0.45\textwidth]{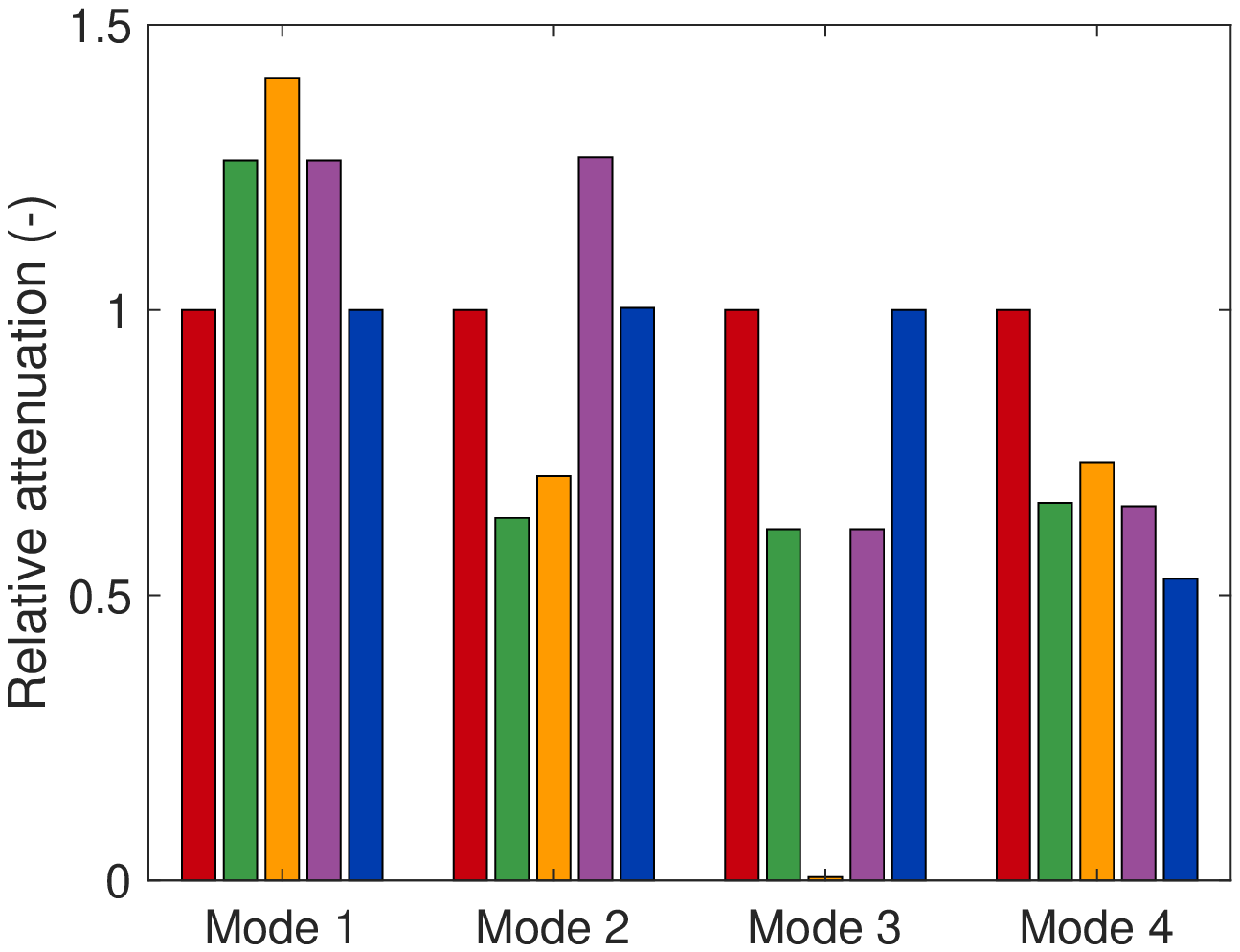}
            \caption{Relative attenuation (compared to the baseline case) on each mode: $[\overline{d}_{p,1}, \overline{d}_{p,2} ,\overline{d}_{p,3} ,\overline{d}_{p,4}]= [1,1,1,1]$ (\textbf{\textcolor{col2}{$\blacksquare$}}), $[\overline{d}_{p,1}, \overline{d}_{p,2} ,\overline{d}_{p,3} ,\overline{d}_{p,4}]= [2,1,1,1]$ (\textbf{\textcolor{col3}{$\blacksquare$}}), $[\overline{d}_{p,1}, \overline{d}_{p,2} ,\overline{d}_{p,3} ,\overline{d}_{p,4}]= [2,1,0,1]$ (\textbf{\textcolor{col5}{$\blacksquare$}}), $[\overline{d}_{p,1}, \overline{d}_{p,2} ,\overline{d}_{p,3} ,\overline{d}_{p,4}]= [2,2,1,1]$ (\textbf{\textcolor{col4}{$\blacksquare$}}) and $[\overline{d}_{p,1}, \overline{d}_{p,2} ,\overline{d}_{p,3} ,\overline{d}_{p,4}]= [2,2,2,1]$ (\textbf{\textcolor{col1}{$\blacksquare$}}).}
            \label{fig:beamAttenuation}
        \end{figure}

        The influence of the relative scaling factors can also be studied. To make things easier to understand, the piezoelectric structure is simplified by grouping the 20 patches into two groups of ten adjacent patches and connecting the groups in parallel. A first network can be synthesized with unit relative scaling factors ($\overline{d}_{p,k}=1$, $k=1,2,3,4$). The resulting FRF is shown in \autoref{fig:beamFRFScaling}. Compared to \autoref{fig:beamFRF0}, the parallel connection of the patches has reduced the control authority of the network, especially on mode 4. This can be understood by the fact that the parallel connection of the patches makes this mode almost unobservable and uncontrollable because of charge cancellation. This baseline case is compared in \autoref{fig:beamAttenuation} to every other investigated case studied hereafter in terms of attenuation on each mode.
        
        In an attempt to improve the control authority over mode 1, one may set $\overline{d}_{p,1} = 2$ while leaving the other relative scaling factors unchanged. In doing so, \autoref{fig:beamFRFScaling} shows that the vibration reduction on mode 1 can be improved by 2dB, but this is done at the expense of vibration reduction on the other modes (by approximately 4dB for all of them). Thus, the relative scaling factors can be used to favor vibration mitigation over some modes, as expected. The vibration reduction improvement on mode 1 can be maximized if mode 3 is left uncontrolled. In this case, a further 1dB can be obtained on mode 1.
        
\begin{figure}[!ht]
            \centering
			\begin{subfigure}[c]{0.5\textwidth}
		         \centering
                \begin{subfigure}[c]{0.48\textwidth}
		            \centering
                    \includegraphics[width=\textwidth]{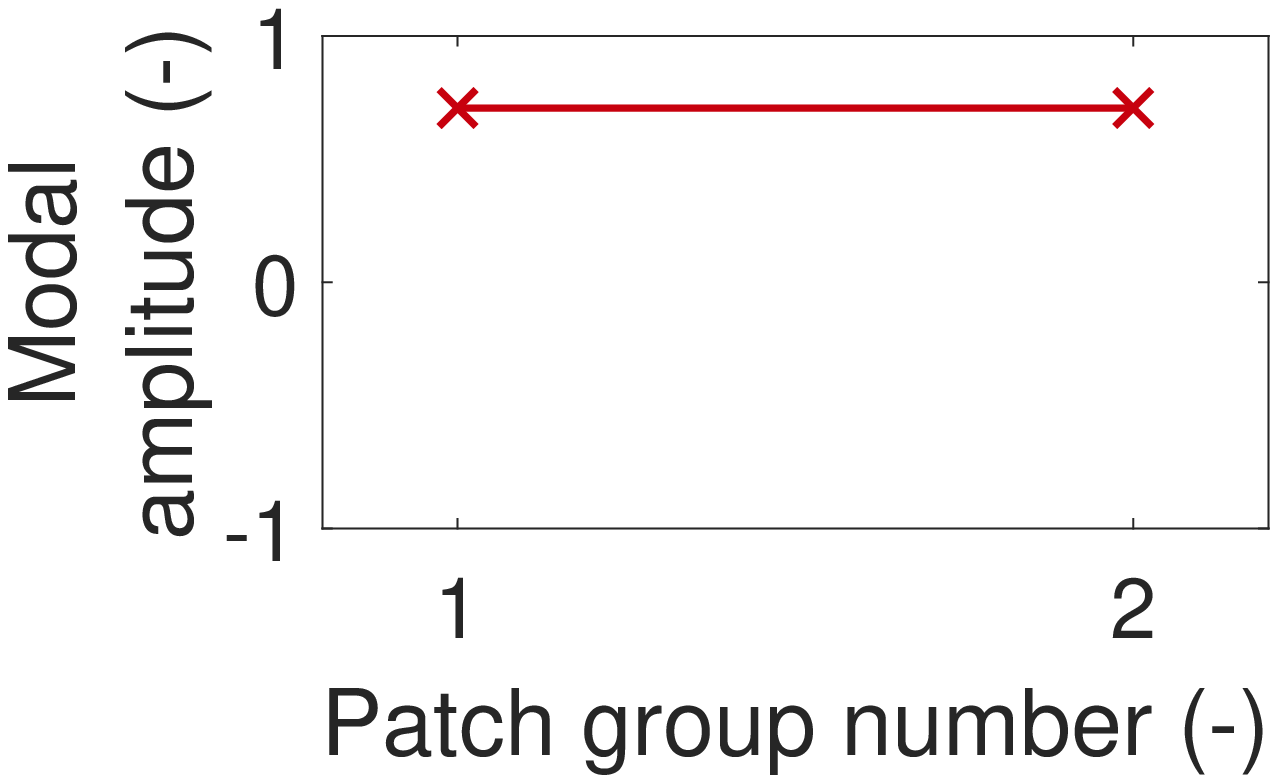}
                    \caption{}
                    \label{sfig:beamScalingMode1}
                \end{subfigure}
                \hspace{0.01\textwidth}
                \begin{subfigure}[c]{0.48\textwidth}
		            \centering
                    \includegraphics[width=\textwidth]{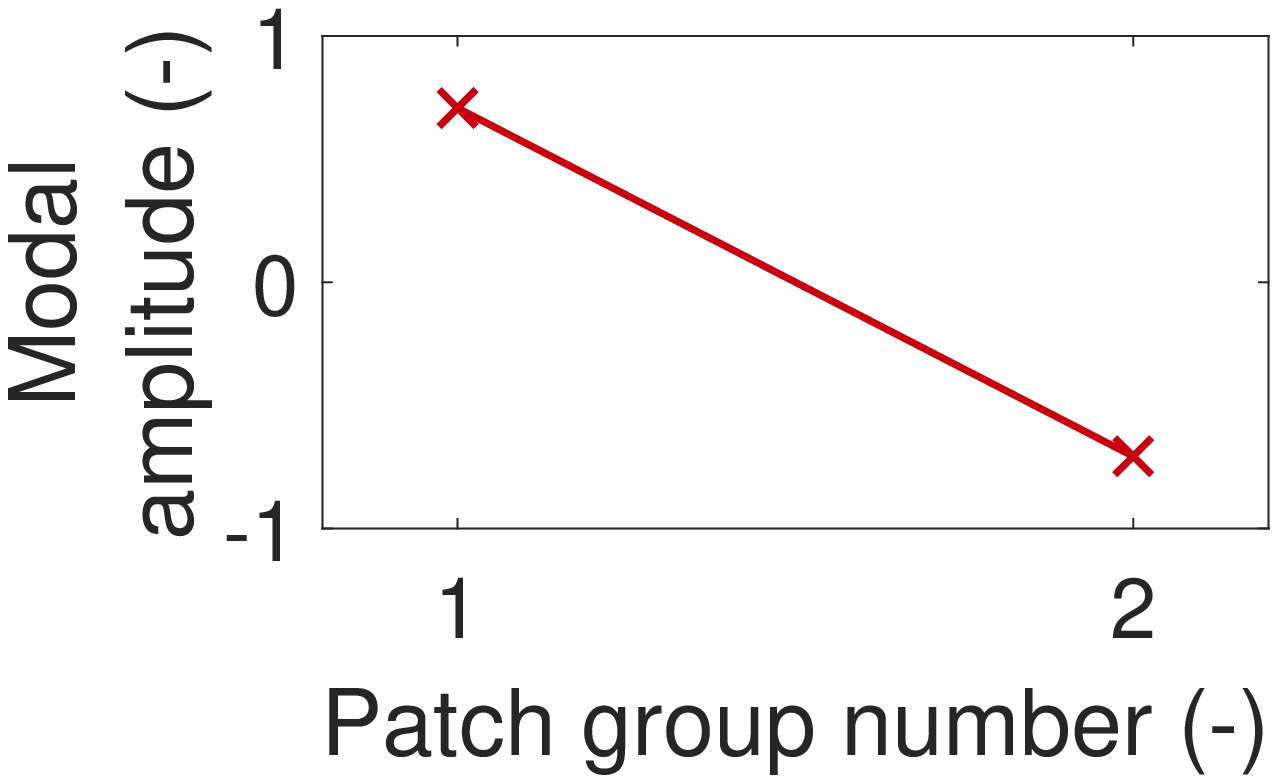}
                    \caption{}
                    \label{sfig:beamScalingMode2}
                \end{subfigure}
                \begin{subfigure}[c]{0.48\textwidth}
		            \centering
                    \includegraphics[width=\textwidth]{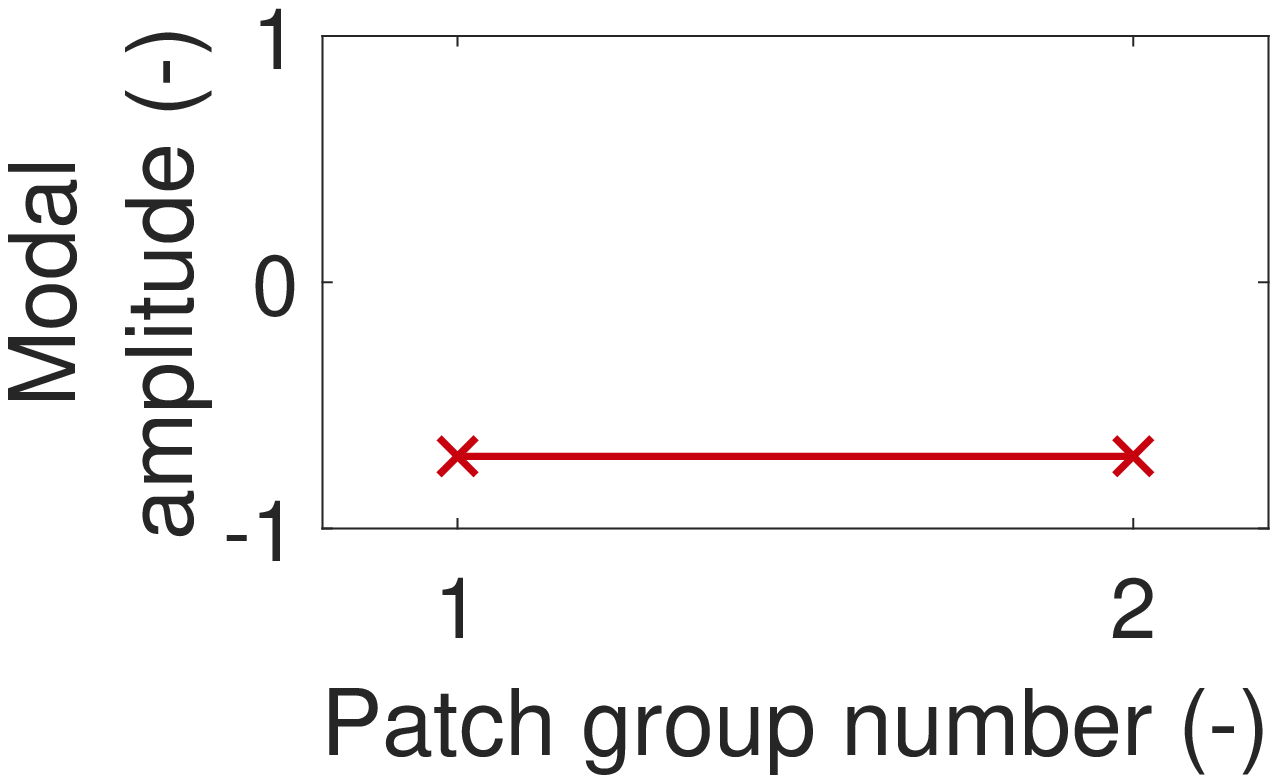}
                    \caption{}
                    \label{sfig:beamScalingMode3}
                \end{subfigure}
                \hspace{0.01\textwidth}
                \begin{subfigure}[c]{0.48\textwidth}
		            \centering
                    \includegraphics[width=\textwidth]{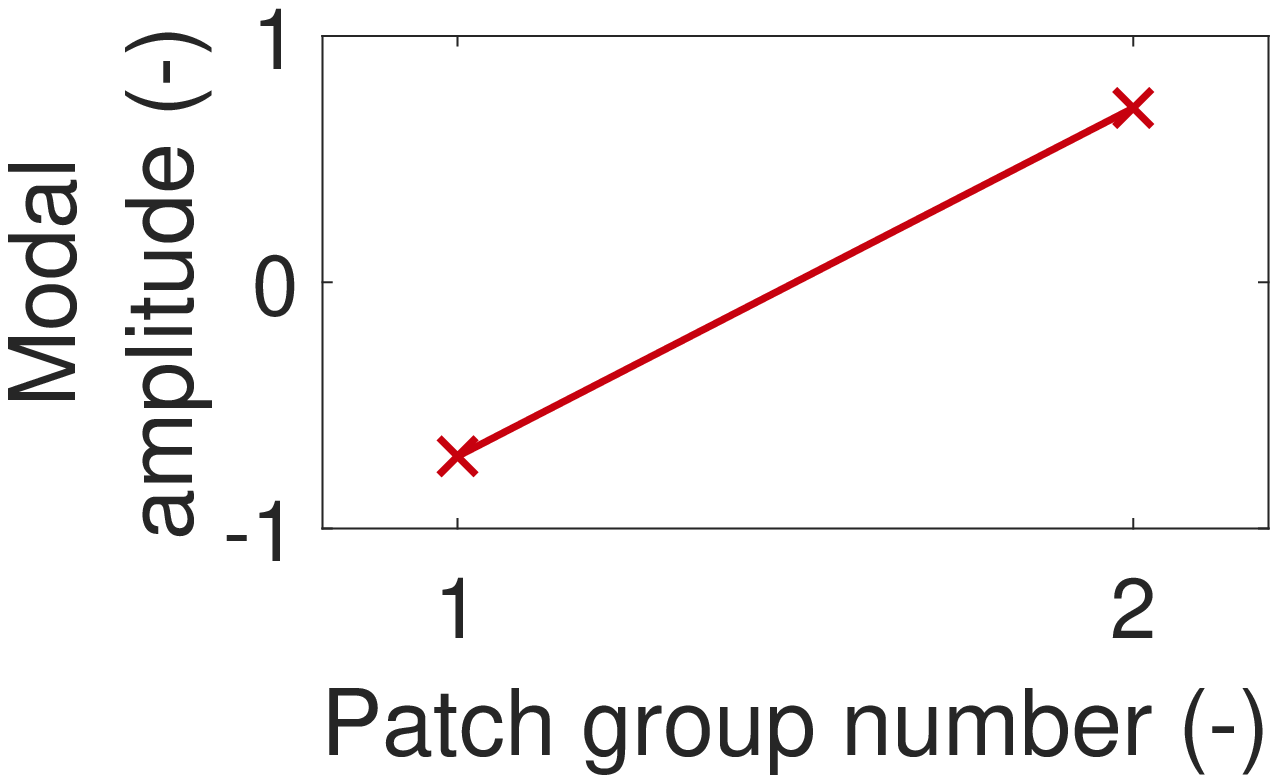}
                    \caption{}
                    \label{sfig:beamScalingMode4}
                \end{subfigure}
            \end{subfigure}
            \hspace{0.02\textwidth}
			\begin{subfigure}[c]{0.45\textwidth}
		        \centering
            \includegraphics[width=\textwidth]{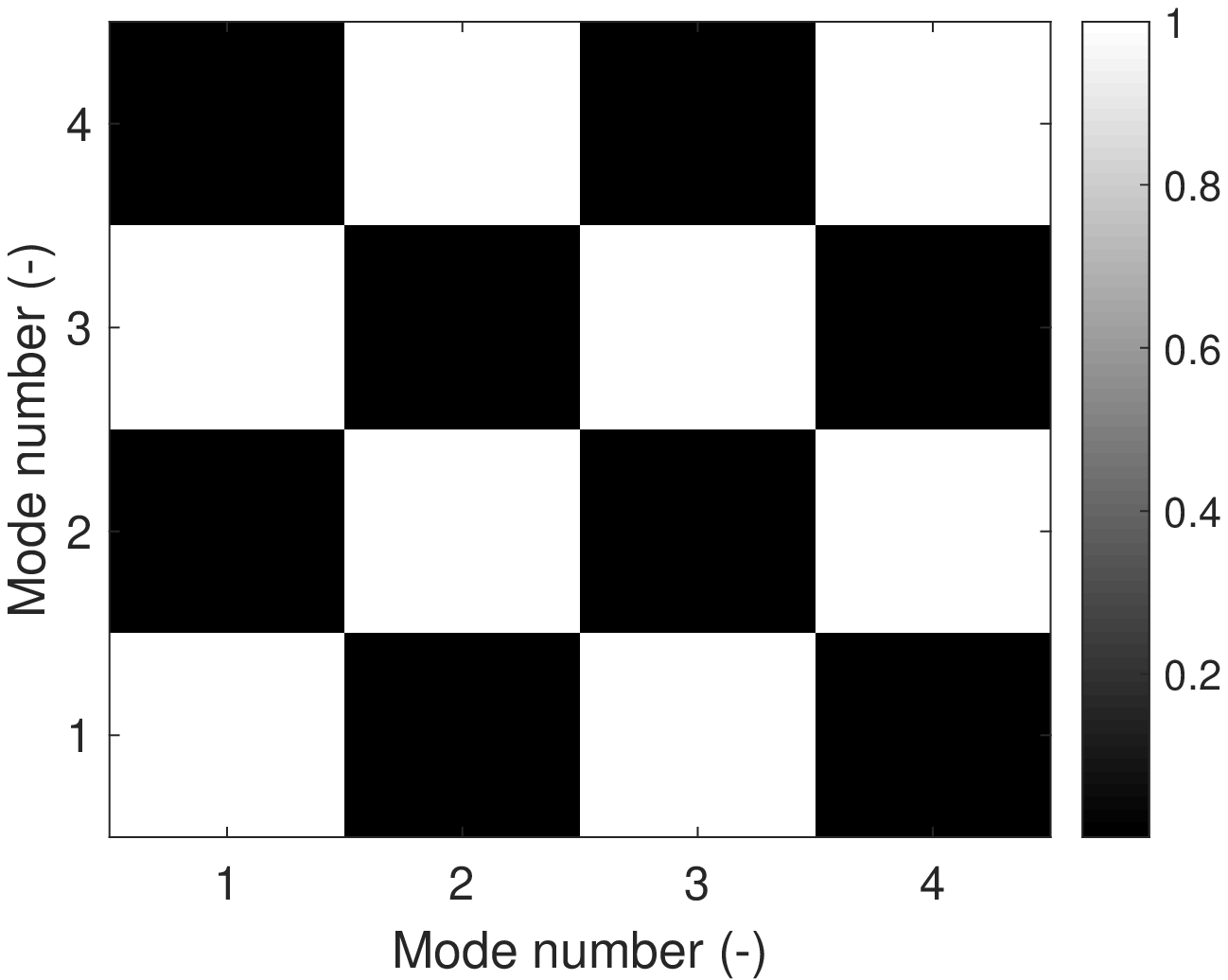}
            \caption{}
            \label{sfig:beamAutoMac}
            \end{subfigure}
            \caption{Optimal dimensionless mode shape 1~\subref{sfig:beamScalingMode1}, 2~\subref{sfig:beamScalingMode2}, 3~\subref{sfig:beamScalingMode3}, 4~\subref{sfig:beamScalingMode4} and auto MAC of the optimal dimensionless electrical mode shapes~\subref{sfig:beamAutoMac}.}
            \label{fig:beamScalingModes}
        \end{figure}

        In addition to improving performance on mode 1, mode 2 can also be emphasized by choosing $\overline{d}_{p,1} = \overline{d}_{p,2} = 2$ while leaving the other relative scaling factors equal to one. In this case, \autoref{fig:beamFRFScaling} indicates that the vibration reduction on mode 2 can be improved without affecting the performance on other modes. While this might seem in contradiction with the previous observation, this result can be understood by looking at the modal assurance criterion of the optimal dimensionless electrical mode shapes $\left|\left(\overline{\bm{\Upphi}}_p^\star\right)^T\overline{\bm{\Upphi}}_p^\star\right|^2$ (auto MAC), indicating how correlated these modes are~\cite{Darleux2020}. The dimensionless mode shapes and auto MAC matrix are featured in \autoref{fig:beamScalingModes}, which clearly shows that modes are similar by pairs, as viewed by the connected patches: mode 1 is similar to mode 3 and mode 2 is similar to mode 4. These mode pairs are also orthogonal to each other. Thus, the scaling factors $\overline{d}_{p,1}$ and $\overline{d}_{p,2}$ affect orthogonal directions and thus different eigenvalues of $\overline{\bm{\Upphi}}_p^\star \overline{\mathbf{D}}_p^2\left(\overline{\bm{\Upphi}}_p^\star\right)^T$. When $\overline{d}_{p,2}$ was equal to one, the largest eigenvalue was the one associated to $\overline{d}_{p,1}$. Increasing $\overline{d}_{p,2}$ up to two does not change this fact, and according to \autoref{eq:alphaScaling} $\alpha$ stays identical, and so do $d_{p,1}$, $d_{p,3}$ and $d_{p,4}$, while $d_{p,2}$ increases (because $\overline{d}_{p,2}$ increases). Thus, performance on mode 2 (or mode 4) can be improved without affecting mode 1, 3 and 4 (or mode 2) up to some point.

        A final example where mode 3 is also emphasized is shown in \autoref{fig:beamFRFScaling}, where $\overline{d}_{p,1}=\overline{d}_{p,2}=\overline{d}_{p,3}=2$ and $\overline{d}_{p,4}=1$. In this case, performance of mode 1 is affected as expected, because the dimensionless electrical modes 1 and 3 are not orthogonal. Performance on modes 2 and 4 is also affected because the maximum eigenvalue of $\overline{\bm{\Upphi}}_p^\star \overline{\mathbf{D}}_p^2\left(\overline{\bm{\Upphi}}_p^\star\right)^T$ has increased and thus $\alpha$ has decreased. Finally, setting all the relative scaling factors equal to 2 would be identical to setting them to 1. 
        
        This analysis showed that the relative scaling factors can indeed be used to balance the control authority over the modes, but they have to be assigned with care. The most significant scaling factor is the one which affects the largest eigenvalue of $\overline{\bm{\Upphi}}_p^\star \overline{\mathbf{D}}_p^2\left(\overline{\bm{\Upphi}}_p^\star\right)^T$. Orthogonality of the dimensionless electrical mode shapes can be assessed with the auto MAC in order to guide the choice of these scaling factors and to understand the trends in the resulting FRFs. In general, starting by setting them all equal to one would be advisable, and fine-tuning of these factors is possible if the analysis is made cautiously.
        
        \subsubsection{Comparison with an analog electrical network}
        
\begin{figure}[!ht]
            \centering
            \begin{subfigure}[c]{\textwidth}
		         \centering
		         \includegraphics[width=\textwidth]{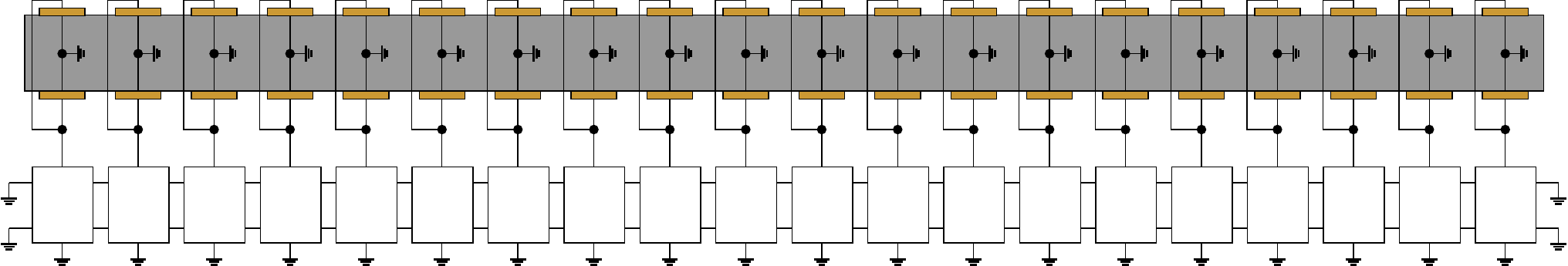}
		         \caption{}
		         \label{sfig:beam}
		     \end{subfigure}
		     \begin{subfigure}[c]{0.45\textwidth}
		         \centering
		         \includegraphics[width=\textwidth]{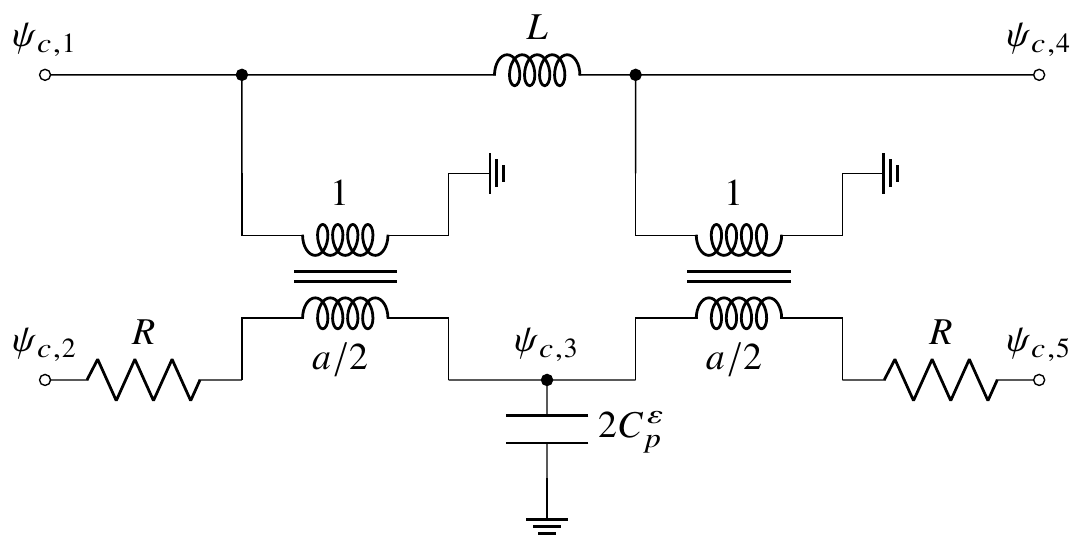}
		         \caption{}
		         \label{sfig:beamAnalog}
		     \end{subfigure}
            \caption{Schematic representation of a free-free beam (in gray) coupled to an electrical network (in white) through piezoelectric patches (in orange)~\subref{sfig:beam} and electrical cell analog to the mechanical cell~\subref{sfig:beamAnalog}.}
            \label{fig:beamSystem}
        \end{figure}

        An electrical cell having dynamics analog to those of the mechanical cell was proposed by Porfiri et al~\cite{Porfiri2004}; it is shown in \autoref{sfig:beamAnalog}. The electrical matrices of this cell can be built as
        \begin{equation}
            \mathbf{C}_c = \begin{bmatrix} 0 \\ 0 \\ 1 \\ 0 \\ 0 \end{bmatrix} 2C_p^\varepsilon \begin{bmatrix} 0 \\ 0 \\ 1 \\ 0 \\ 0 \end{bmatrix}^T, \quad 
            \mathbf{G}_c = \begin{bmatrix} -a \\ 1 \\ -1 \\ 0 \\ 0 \end{bmatrix}\dfrac{1}{R}\begin{bmatrix} -a \\ 1 \\ -1 \\ 0 \\ 0 \end{bmatrix}^T + \begin{bmatrix} 0 \\ 0 \\ -1 \\ a \\ 1 \end{bmatrix}\dfrac{1}{R}\begin{bmatrix} 0 \\ 0 \\ -1 \\ a \\ 1 \end{bmatrix}^T , \quad
            \mathbf{B}_c = \begin{bmatrix} 1 \\ 0 \\ 0 \\ -1 \\ 0 \end{bmatrix} \dfrac{1}{L} \begin{bmatrix} 1 \\ 0 \\ 0 \\ -1 \\ 0 \end{bmatrix}^T,
        \end{equation}
        using the same ordering in the degrees of freedom as in \autoref{sfig:beamAnalog}~\cite{Grainger1994}. The localization matrix in this cell is $\mathbf{E}_{p,c}^T=\begin{bmatrix}0 & 0 & 1 & 0 & 0\end{bmatrix}$ in this case. The characteristics of the electrical cell are reported in \autoref{tab:cellCharacs}, where the resistance was optimally tuned to the first mode. The matrices of the overall network $\mathbf{C}$, $\mathbf{G}$, $\mathbf{B}$ and $\mathbf{E}_p$ can then be built by standard finite element assembly procedures~\cite{Geradin2014}. The analog network features the same resonance frequencies and mode shapes as the beam. When the two systems are coupled via the piezoelectric patches, as shown in \autoref{sfig:beam}, broadband damping in the structure is achieved.
        
        \begin{table}[!ht]
        \centering
        \caption{Characteristics of the electrical cell.}
        \label{tab:cellCharacs}
        \begin{tabular}{ccccc}
        \hline
        Parameter & $C_p^\varepsilon$ & $R$ & $L$ & $a$ \\
        \hline
        Value & 21.96nF & 57.5$\Upomega$ & 161.1mH & 1 \\
        \hline
        \end{tabular}
        \vspace*{-4pt}
        \end{table}

        \begin{figure}[!ht]
            \centering
            \includegraphics[width=0.8\textwidth]{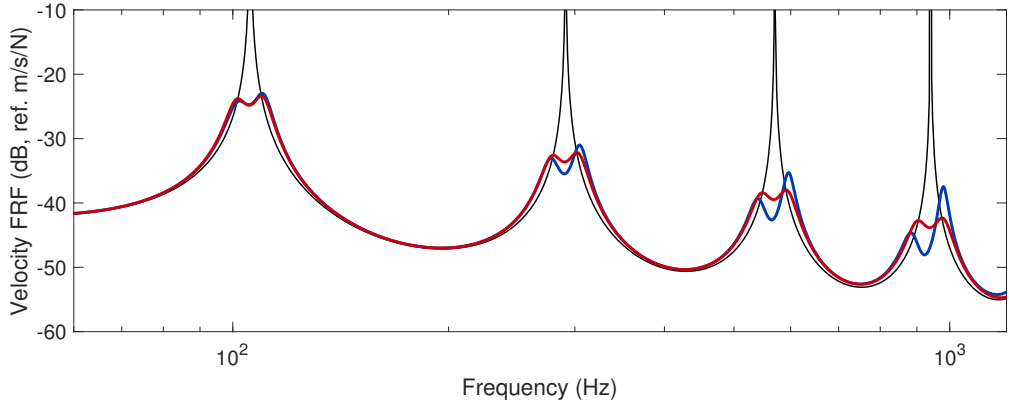}
            \caption{Velocity FRF of the beam with short-circuited patches~(---), with a network assembled from analog electrical cells~\cite{Lossouarn2015a}~(\textbf{\textcolor{col1}{---}}) and with a network synthesized with the modal-based approach~(\textbf{\textcolor{col2}{---}}).}
            \label{fig:beamFRF}
        \end{figure}
        
         \autoref{fig:beamFRF} indicates that the modal-based approach yields a more accurate tuning of the electrical network with respect to the $H_\infty$ norm. However, unlike the analog network based on elementary cells, the electrical matrices do not have a band structure, which means that there may exist a large number of interconnections in the network. Thus, the network obtained with the proposed approach may potentially be more difficult to realize practically. 
        
\begin{figure}[!ht]
            \centering
			\begin{subfigure}[c]{0.45\textwidth}
		         \centering
                \includegraphics[width=\textwidth]{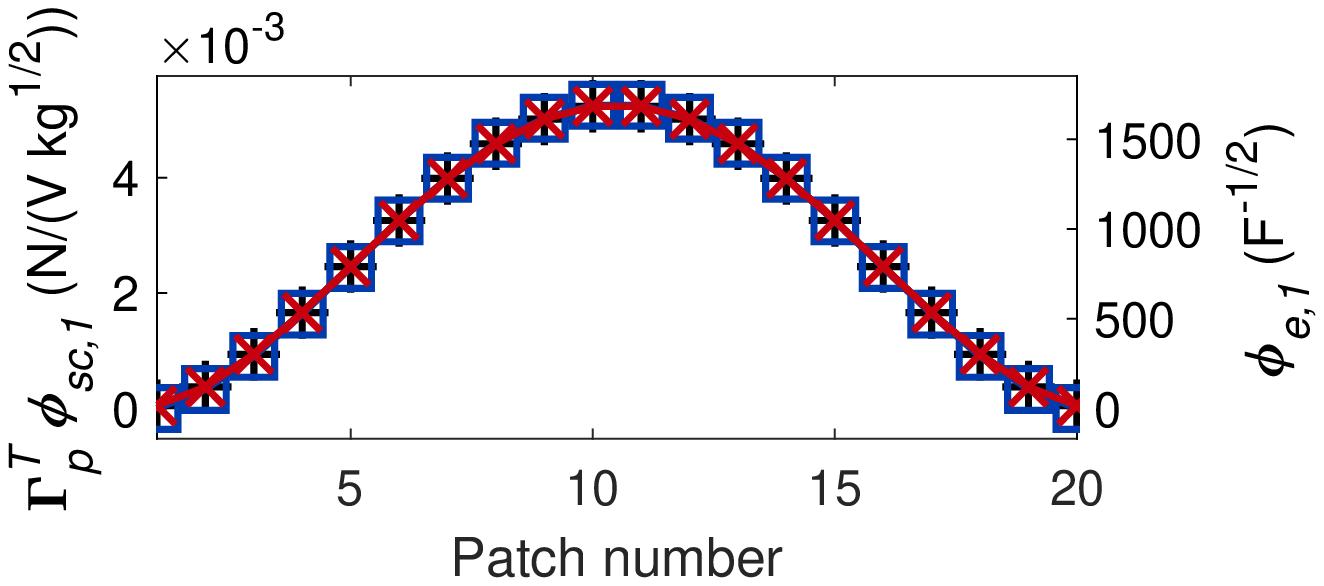}
                \caption{}
                \label{sfig:beamMode1}
            \end{subfigure}
            \hspace{0.05\textwidth}
			\begin{subfigure}[c]{0.45\textwidth}
		         \centering
                \includegraphics[width=\textwidth]{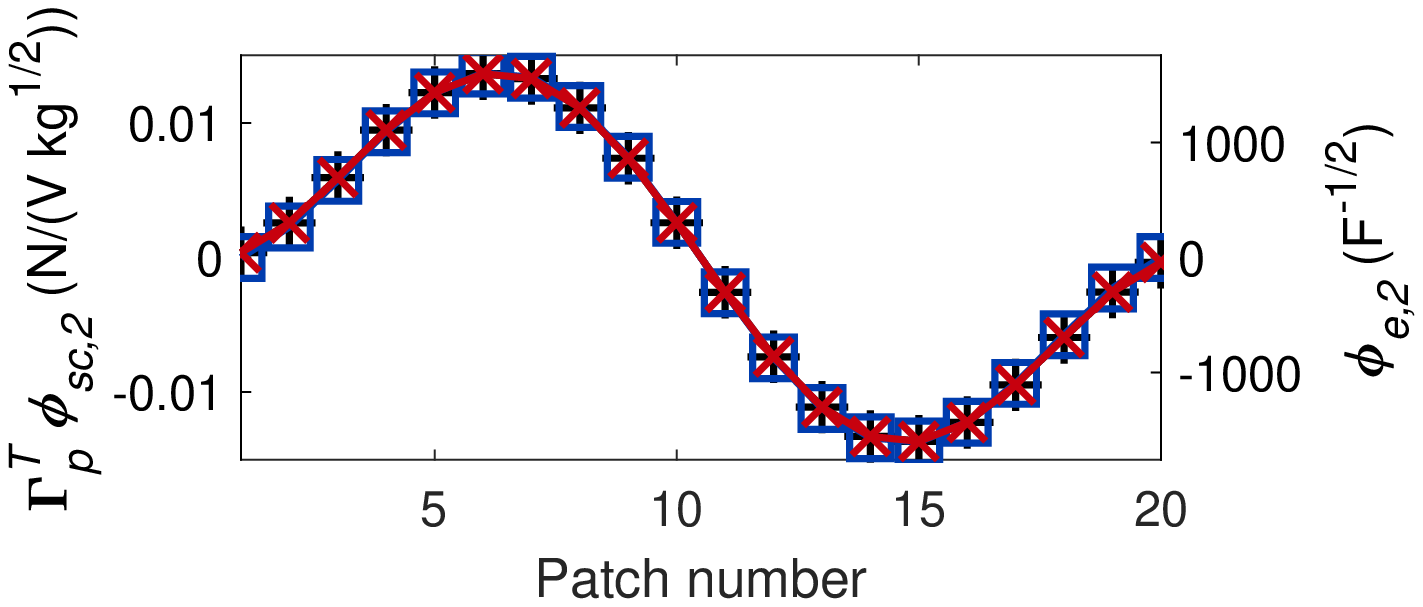}
                \caption{}
                \label{sfig:beamMode2}
            \end{subfigure}
			\begin{subfigure}[c]{0.45\textwidth}
		         \centering
                \includegraphics[width=\textwidth]{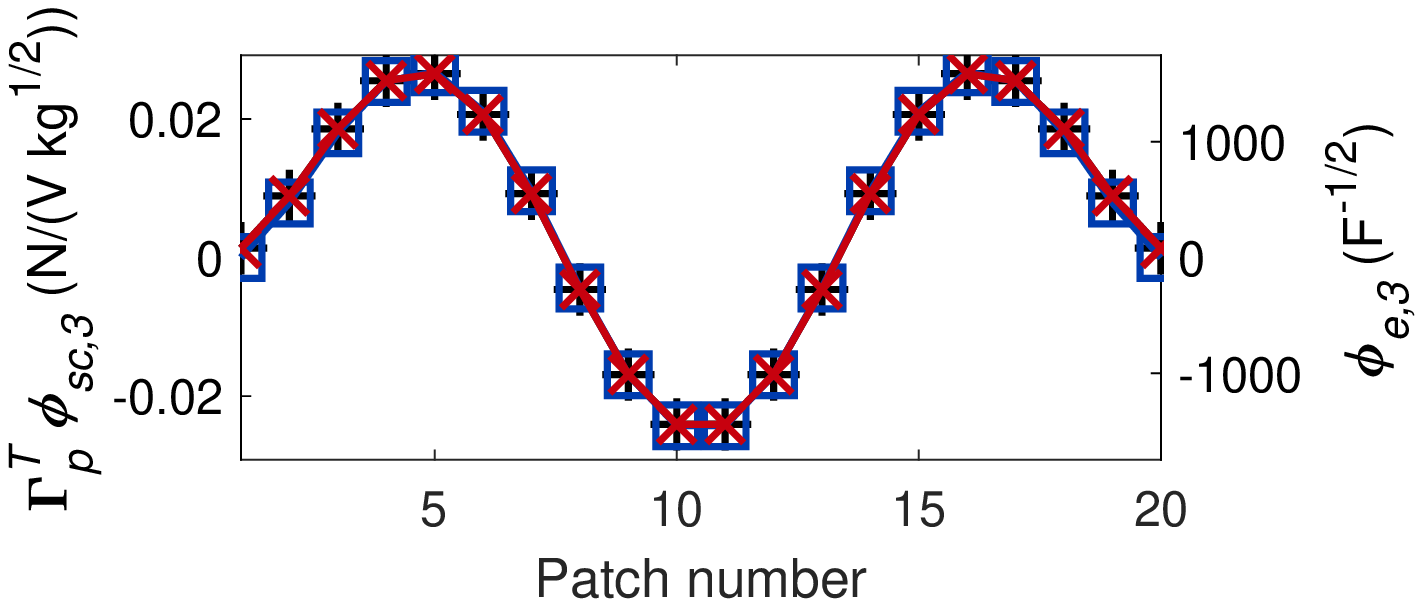}
                \caption{}
                \label{sfig:beamMode3}
            \end{subfigure}
            \hspace{0.05\textwidth}
			\begin{subfigure}[c]{0.45\textwidth}
		         \centering
                \includegraphics[width=\textwidth]{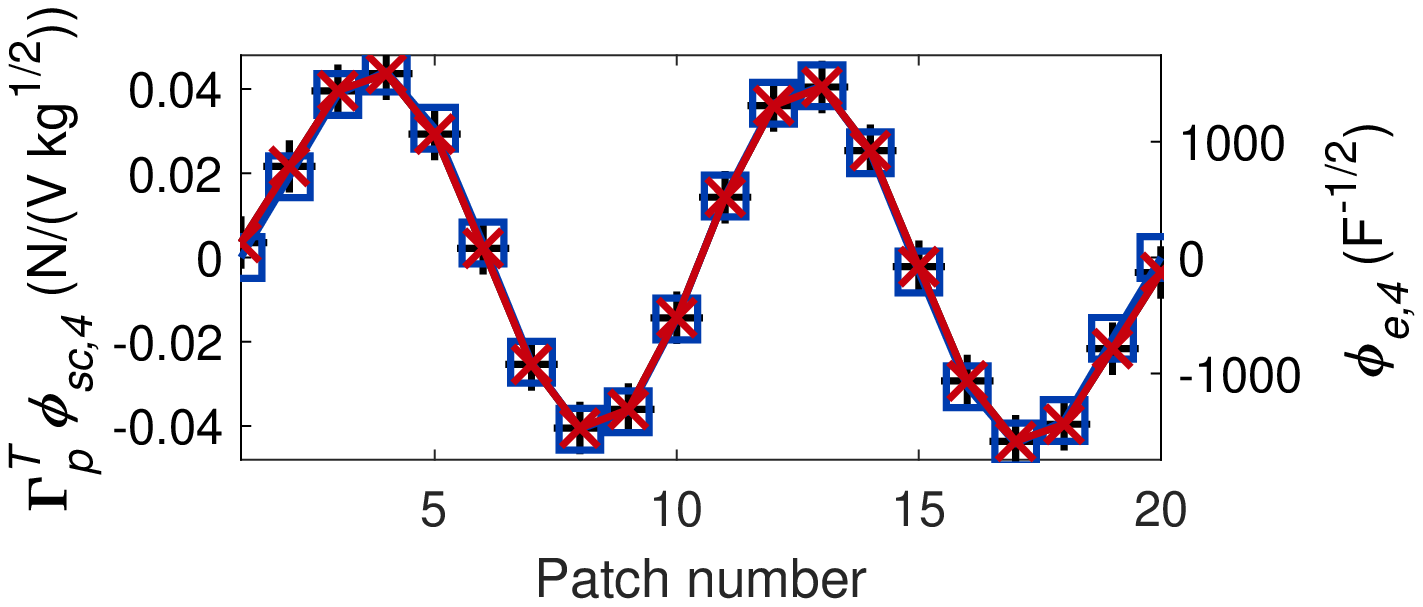}
                \caption{}
                \label{sfig:beamMode4}
            \end{subfigure}
            \caption{Modal strain in the piezoelectric transducers (\textbf{-+-}), electrical mode shapes of a network assembled from analog electrical cells~\cite{Lossouarn2015a} (\textbf{\textcolor{col1}{-$\square$-}}) and electrical mode shapes of a network synthesized with the modal-based approach (\textbf{\textcolor{col2}{-{\scriptsize $\times$}-}}): mode 1~\subref{sfig:beamMode1}, 2~\subref{sfig:beamMode2}, 3~\subref{sfig:beamMode3} and 4~\subref{sfig:beamMode4}.}
            \label{fig:beamModes}
        \end{figure}

        It should be pointed out that since the piezoelectric patches are all identical, the network synthesized with the modal-based approach is also an analog of the structure in a prescribed frequency band, because it features identical resonance frequencies and mode shapes as the targeted ones. This is confirmed in \autoref{fig:beamModes}, where the first four piezoelectric modal strain shapes and the first four electrical mode shapes obtained with either method are displayed (the patch are numbered according to their position, from one end of the beam to the other).
        
    \subsection{Fully clamped plate}
        
\begin{figure}[!ht]
            \centering
			\begin{subfigure}[c]{0.6\textwidth}
		         \centering
                    \includegraphics[width=\textwidth]{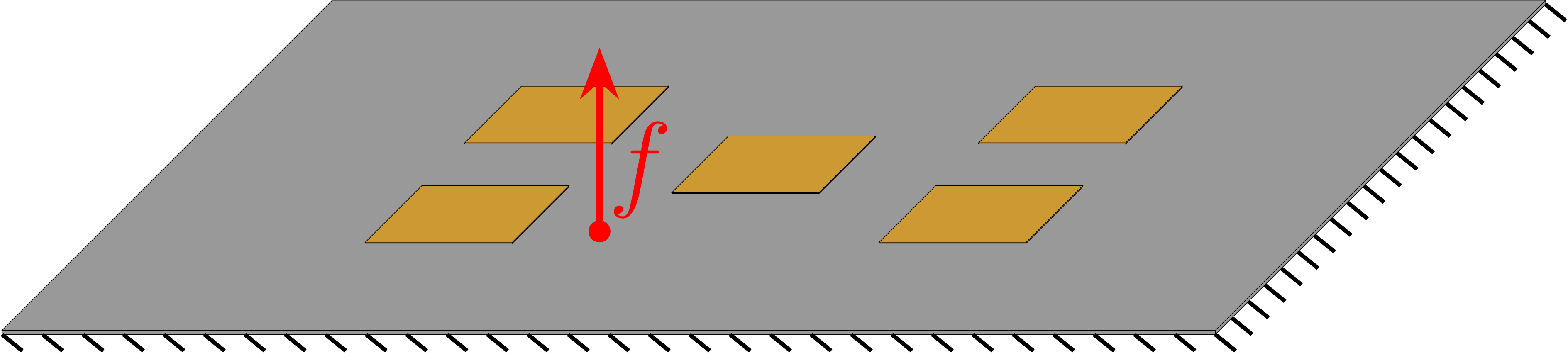}
		         \caption{}
		         \label{sfig:plate_schematics}
		     \end{subfigure}
		     \hspace{0.05\textwidth}
			\begin{subfigure}[c]{0.3\textwidth}
		         \centering
                    \includegraphics[width=\textwidth]{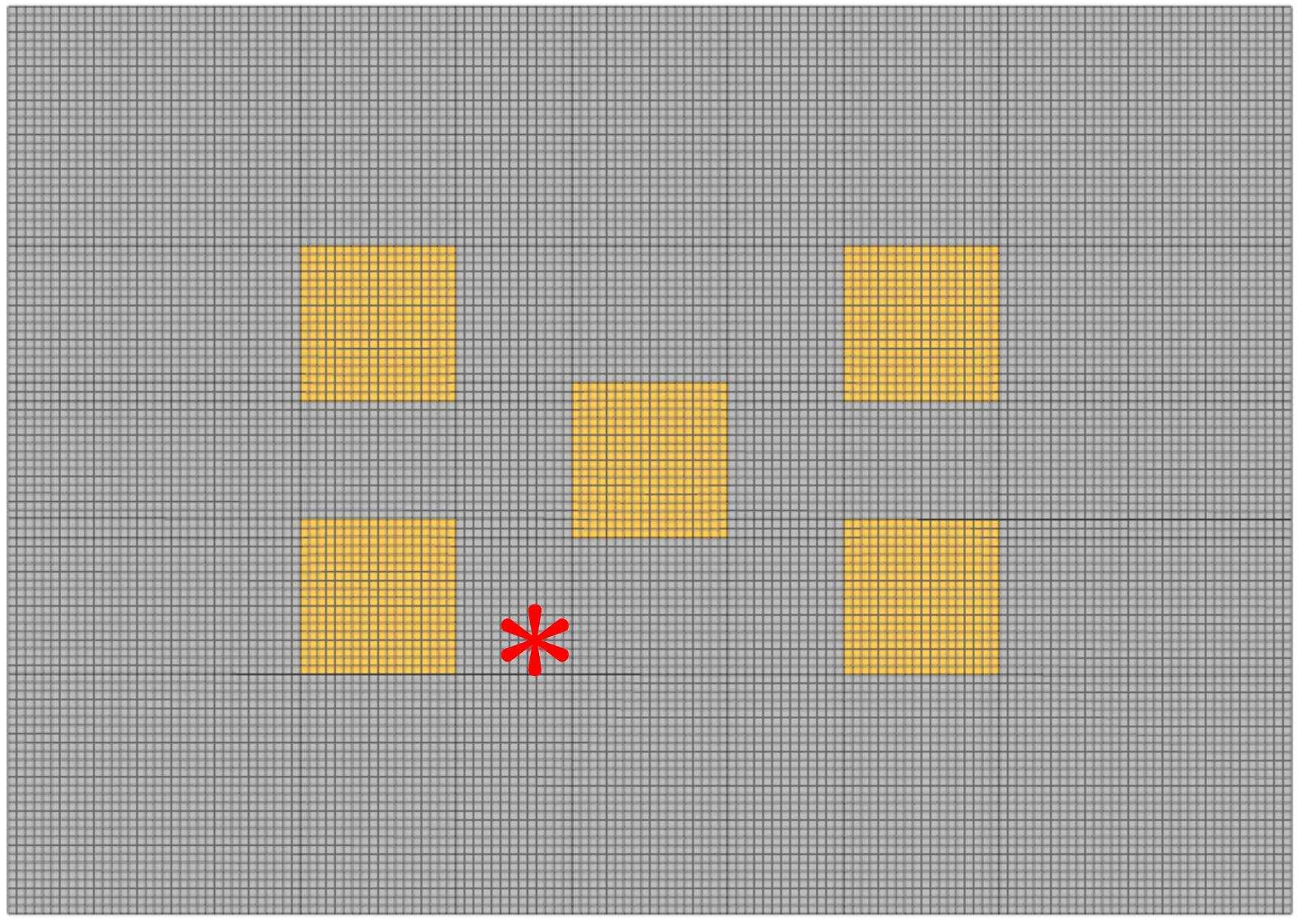}
		         \caption{}
		         \label{sfig:plate_FEM}
		     \end{subfigure}
            \caption{Schematic representation of a fully clamped plate (in gray) excited by a point force (in red) to which are bonded piezoelectric patches (in orange) \subref{sfig:plate_schematics} and finite element mesh of the plate \subref{sfig:plate_FEM}.}
            \label{fig:plate}
        \end{figure}

        The second example is the fully clamped plate depicted in \autoref{sfig:plate_schematics}~\cite{Giorgio2009}. Five pairs of piezoelectric patches are bonded symmetrically on either side of the plate, each pair being electrically connected in parallel. The finite element model of the piezoelectric structure in \autoref{sfig:plate_FEM} was developed in SAMCEF (see~\cite{Piefort2001} for implementation details) and a reduced-order model was imported in MATLAB. The Craig-Bampton technique was used with 20 retained modes~\cite{Geradin2014}. The plate is subjected to a pointwise force located on a node of the finite element mesh at 41\% of the plate length and 30\% of the plate width from the lower left corner in \autoref{fig:plate}, in order to make it as close as possible to what was used in~\cite{Giorgio2009}. Besides that, the characteristics of the system were identical to those reported therein.
        
\begin{figure}[!ht]
            \centering
            \includegraphics[width=0.8\textwidth]{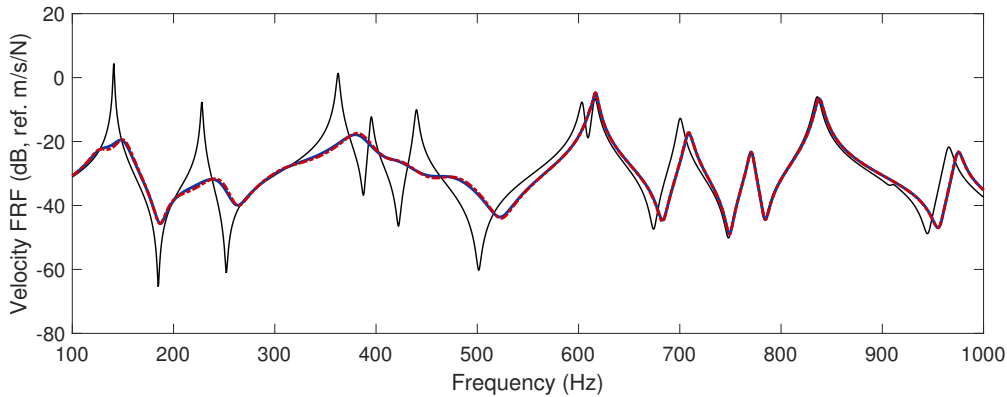}
            \caption{Velocity FRF of the plate with short-circuited patches~(---), with a network synthesized with the method in~\cite{Giorgio2009}~(\textbf{\textcolor{col1}{---}}) and with a network synthesized with the modal-based approach~(\textbf{\textcolor{col2}{-$\cdot$-}}).}
            \label{fig:plateFRF}
        \end{figure}

        The driving-point FRF of the plate coupled to the electrical network obtained with a modal-based synthesis targeting the first five modes of the plate is shown in \autoref{fig:plateFRF}. Identical scaling factors for the five modes were used. As for the beam, the electrical network can very effectively damp out the resonant vibrations of the targeted modes.
        
        The method in~\cite{Giorgio2009} was also used to synthesize an electrical network. The cornerstone of that method is to find an orthogonal transformation matrix $\mathbf{U}$ for the $P=N_s$ electrical degrees of freedom such that $\bm{\uppsi} = \mathbf{U}\bm{\upchi}$ (where $\bm{\upchi}$ are the transformed electrical degrees of freedom). The matrix $\mathbf{U}$ is aimed to make the transformed piezoelectric coupling matrix as close as possible from a diagonal matrix. The unknowns of this problem are the $N_s^2$ entries of $\mathbf{U}$. A system of $N_s^2$ quadratic equations is built to impose the orthogonality of $\mathbf{U}$ and the optimal closeness of the transformed coupling matrix to a diagonal matrix. It can be solved numerically to find the unknowns. Similarly to that method, the present paper uses a transformation of the electrical degrees of freedom with the electrical modes, $\bm{\uppsi} = \bm{\Upphi}_e\bm{\eta}_e$, but the matrix $\bm{\Upphi}_e$ needs not be orthogonal. It is not explicitly set to diagonalize the piezoelectric coupling matrix, and it can be computed directly. Moreover, the modal-based method does not require the number of controlled modes $N_s$ to be equal to the number of transducers $P$. \autoref{fig:plateFRF} reveals that the performance of the resulting network is found to be almost identical to that of our approach. Despite the marked implementation differences between the methods, they both aim at optimally using the control capability offered by the transducers, which makes this result quite expectable in the end.
        
\begin{figure}[!ht]
            \centering
            \includegraphics[width=0.8\textwidth]{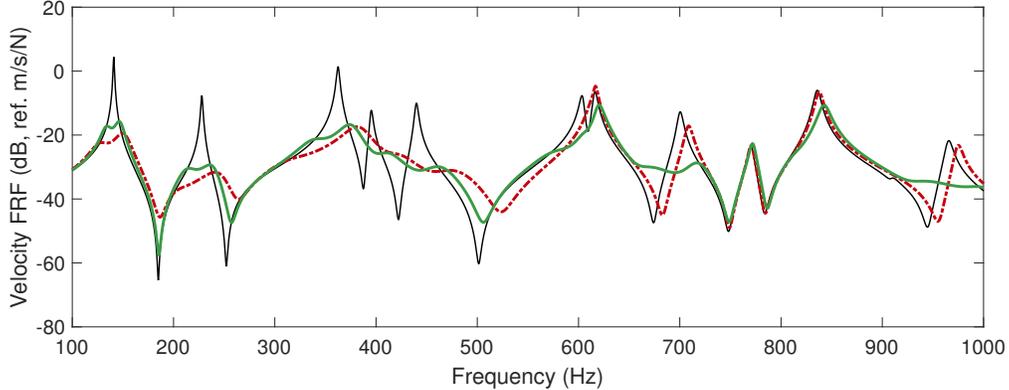}
            \caption{Velocity FRF of the plate with short-circuited patches~(---) and with a network synthesized with the modal-based approach targeting the first five~(\textbf{\textcolor{col2}{-$\cdot$-}}) and twelve~(\textbf{\textcolor{col3}{---}}) modes.}
            \label{fig:plateFRF2}
        \end{figure}

\begin{figure}[!ht]
            \centering
            \includegraphics[width=0.45\textwidth]{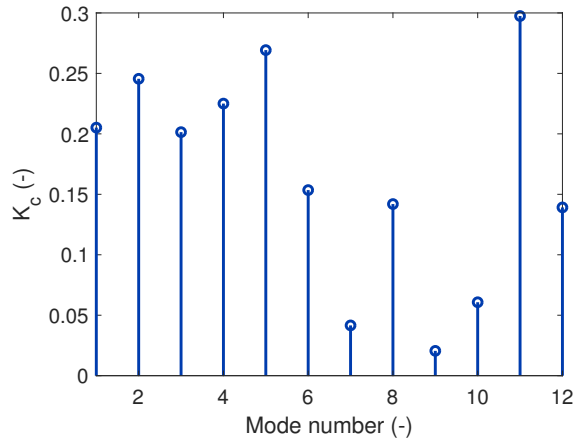}
            \caption{EEMCF of the first fifteen modes of the plate obtained by short- and open-circuiting all the patches.}
            \label{fig:plateKc}
        \end{figure}

        A network targeting more modes than the number of piezoelectric patches was also designed. In order mitigate the vibration between 0 and 1000Hz, the first twelve modes with identical scaling factors were considered. The network performance is compared against that of the network targeting five modes in \autoref{fig:plateFRF2}. While vibration mitigation is less effective for the first few modes, the higher-frequency modes are now damped with the exception of the 7th, 9th and 10th modes. This is explained by the fact that the electromechanical coupling that the patches have with these modes is rather low, as shown in \autoref{fig:plateKc}. In this plot, the EEMCFs were obtained by considering the resonance frequencies of the plate when all the patches are shorted, and when all the patches are in open-circuit. Better mitigation results could be obtained with more patches or with optimally-located patches.
        
\section{Conclusion}
\label{sec:conclusion}

Aiming to provide a simple and systematic design strategy to mitigate multiple resonances of complex structures, this work leveraged the concept of passive electrical networks interconnecting piezoelectric transducers. The core idea is to tailor the modal properties of the electrical network. Specifically, the key original aspect of this study is to guarantee the passivity of the network interconnecting the transducers through mathematical inequalities involving the electrical mode shapes and the capacitance matrix of the transducers at constant strain. The second originality is then to optimize the modal shapes for the maximization of the modal electromechanical coupling coefficients. Eventually, the matrices governing the behavior of the electrical network in the physical domain can be retrieved from the knowledge of the modal properties by simple inversions.

This modal-based synthesis was illustrated using two examples, namely a free-free beam and a fully clamped plate. In both cases, effective broadband damping could be obtained with performance comparable to that of the state-of-the-art approaches. 

This work opens up interesting perspectives. A first one concerns the practical realization of the electrical network, i.e., to find the actual electrical components to implement it. An alternative and attractive solution is to resort to microcontrollers as in, e.g., \cite{Giorgio2009,Raze2020c}, for which virtually any admittance can be implemented. A second perspective would be to study how the design choices (such as the scaling factors and electrical mode shapes on the internal degrees of freedom) impact the topology, number of components and component values of the network. Finally, this approach could also be experimentally validated for broadband damping of real-world structures, such as, e.g., bladed assemblies.

	\bibliography{bibliography}

\begin{thebibliography}{10}
\providecommand{\url}[1]{#1}
\csname url@samestyle\endcsname
\providecommand{\newblock}{\relax}
\providecommand{\bibinfo}[2]{#2}
\providecommand{\BIBentrySTDinterwordspacing}{\spaceskip=0pt\relax}
\providecommand{\BIBentryALTinterwordstretchfactor}{4}
\providecommand{\BIBentryALTinterwordspacing}{\spaceskip=\fontdimen2\font plus
\BIBentryALTinterwordstretchfactor\fontdimen3\font minus
  \fontdimen4\font\relax}
\providecommand{\BIBforeignlanguage}[2]{{%
\expandafter\ifx\csname l@#1\endcsname\relax
\typeout{** WARNING: IEEEtran.bst: No hyphenation pattern has been}%
\typeout{** loaded for the language `#1'. Using the pattern for}%
\typeout{** the default language instead.}%
\else
\language=\csname l@#1\endcsname
\fi
#2}}
\providecommand{\BIBdecl}{\relax}
\BIBdecl

\bibitem{Hagood1991}
\BIBentryALTinterwordspacing
N.~Hagood and A.~von Flotow, ``{Damping of structural vibrations with
  piezoelectric materials and passive electrical networks},'' \emph{Journal of
  Sound and Vibration}, vol. 146, no.~2, pp. 243--268, apr 1991. [Online].
  Available:
  \url{https://linkinghub.elsevier.com/retrieve/pii/0022460X91907629}
\BIBentrySTDinterwordspacing

\bibitem{Wu1996}
\BIBentryALTinterwordspacing
S.-y. Wu, ``{Piezoelectric shunts with a parallel R-L circuit for structural
  damping and vibration control},'' in \emph{Spie}, C.~D. Johnson, Ed., vol.
  2720, may 1996, pp. 259--269. [Online]. Available:
  \url{http://proceedings.spiedigitallibrary.org/proceeding.aspx?articleid=1017503}
\BIBentrySTDinterwordspacing

\bibitem{Moheimani2006}
\BIBentryALTinterwordspacing
S.~O.~R. Moheimani and A.~J. Fleming, \emph{{Piezoelectric Transducers for
  Vibration Control and Damping}}, ser. Advances in Industrial Control.\hskip
  1em plus 0.5em minus 0.4em\relax London: Springer-Verlag, 2006. [Online].
  Available: \url{http://link.springer.com/10.1007/1-84628-332-9}
\BIBentrySTDinterwordspacing

\bibitem{Berardengo2017}
\BIBentryALTinterwordspacing
M.~Berardengo, S.~Manzoni, and A.~M. Conti, ``{Multi-mode passive piezoelectric
  shunt damping by means of matrix inequalities},'' \emph{Journal of Sound and
  Vibration}, vol. 405, pp. 287--305, 2017. [Online]. Available:
  \url{http://dx.doi.org/10.1016/j.jsv.2017.06.002}
\BIBentrySTDinterwordspacing

\bibitem{Raze2020}
\BIBentryALTinterwordspacing
G.~Raze, A.~Paknejad, G.~Zhao, C.~Collette, and G.~Kerschen, ``{Multimodal
  vibration damping using a simplified current blocking shunt circuit},''
  \emph{Journal of Intelligent Material Systems and Structures}, vol.~31,
  no.~14, pp. 1731--1747, aug 2020. [Online]. Available:
  \url{http://journals.sagepub.com/doi/10.1177/1045389X20930103}
\BIBentrySTDinterwordspacing

\bibitem{DellIsola1998}
\BIBentryALTinterwordspacing
F.~Dell'Isola and S.~Vidoli, ``{Continuum modelling of piezoelectromechanical
  truss beams: an application to vibration damping},'' \emph{Archive of Applied
  Mechanics (Ingenieur Archiv)}, vol.~68, no.~1, pp. 1--19, feb 1998. [Online].
  Available: \url{http://link.springer.com/10.1007/s004190050142}
\BIBentrySTDinterwordspacing

\bibitem{Vidoli2000}
\BIBentryALTinterwordspacing
S.~Vidoli and F.~Dell'Isola, ``{Modal coupling in one-dimensional
  electromechanical structured continua},'' \emph{Acta Mechanica}, vol. 141,
  no. 1-2, pp. 37--50, mar 2000. [Online]. Available:
  \url{http://link.springer.com/10.1007/BF01176806}
\BIBentrySTDinterwordspacing

\bibitem{Alessandroni2002}
\BIBentryALTinterwordspacing
S.~Alessandroni, F.~Dell'Isola, and M.~Porfiri, ``{A revival of electric
  analogs for vibrating mechanical systems aimed to their efficient control by
  PZT actuators},'' \emph{International Journal of Solids and Structures},
  vol.~39, no.~20, pp. 5295--5324, oct 2002. [Online]. Available:
  \url{https://linkinghub.elsevier.com/retrieve/pii/S002076830200402X}
\BIBentrySTDinterwordspacing

\bibitem{Maurini2004}
\BIBentryALTinterwordspacing
C.~Maurini, F.~Dell'Isola, and D.~{Del Vescovo}, ``{Comparison of
  piezoelectronic networks acting as distributed vibration absorbers},''
  \emph{Mechanical Systems and Signal Processing}, vol.~18, no.~5, pp.
  1243--1271, sep 2004. [Online]. Available:
  \url{https://linkinghub.elsevier.com/retrieve/pii/S0888327003000827}
\BIBentrySTDinterwordspacing

\bibitem{Porfiri2004}
\BIBentryALTinterwordspacing
M.~Porfiri, F.~Dell'Isola, and F.~M. {Frattale Mascioli}, ``{Circuit analog of
  a beam and its application to multimodal vibration damping, using
  piezoelectric transducers},'' \emph{International Journal of Circuit Theory
  and Applications}, vol.~32, no.~4, pp. 167--198, jul 2004. [Online].
  Available: \url{http://doi.wiley.com/10.1002/cta.273}
\BIBentrySTDinterwordspacing

\bibitem{Lossouarn2015}
\BIBentryALTinterwordspacing
B.~Lossouarn, M.~Aucejo, and J.-F. De{\"{u}}, ``{Multimodal coupling of
  periodic lattices and application to rod vibration damping with a
  piezoelectric network},'' \emph{Smart Materials and Structures}, vol.~24,
  no.~4, p. 045018, apr 2015. [Online]. Available:
  \url{http://dx.doi.org/10.1088/0964-1726/24/4/045018
  https://iopscience.iop.org/article/10.1088/0964-1726/24/4/045018}
\BIBentrySTDinterwordspacing

\bibitem{Lossouarn2015a}
\BIBentryALTinterwordspacing
B.~Lossouarn, J.~F. De{\"{u}}, and M.~Aucejo, ``{Multimodal vibration damping
  of a beam with a periodic array of piezoelectric patches connected to a
  passive electrical network},'' \emph{Smart Materials and Structures},
  vol.~24, no.~11, p. 115037, nov 2015. [Online]. Available:
  \url{https://iopscience.iop.org/article/10.1088/0964-1726/24/11/115037}
\BIBentrySTDinterwordspacing

\bibitem{Lossouarn2016}
\BIBentryALTinterwordspacing
B.~Lossouarn, J.-F. De{\"{u}}, M.~Aucejo, and K.~A. Cunefare, ``{Multimodal
  vibration damping of a plate by piezoelectric coupling to its analogous
  electrical network},'' \emph{Smart Materials and Structures}, vol.~25,
  no.~11, p. 115042, nov 2016. [Online]. Available:
  \url{https://iopscience.iop.org/article/10.1088/0964-1726/25/11/115042}
\BIBentrySTDinterwordspacing

\bibitem{Darleux2020a}
\BIBentryALTinterwordspacing
R.~Darleux, B.~Lossouarn, and J.-f. De{\"{u}}, ``{Broadband vibration damping
  of non-periodic plates by piezoelectric coupling to their electrical
  analogues},'' \emph{Smart Materials and Structures}, vol.~29, no.~5, p.
  054001, may 2020. [Online]. Available:
  \url{https://iopscience.iop.org/article/10.1088/1361-665X/ab7948}
\BIBentrySTDinterwordspacing

\bibitem{Giorgio2015}
\BIBentryALTinterwordspacing
I.~Giorgio, L.~Galantucci, A.~{Della Corte}, and D.~{Del Vescovo},
  ``{Piezo-electromechanical smart materials with distributed arrays of
  piezoelectric transducers: Current and upcoming applications},''
  \emph{International Journal of Applied Electromagnetics and Mechanics},
  vol.~47, no.~4, pp. 1051--1084, jun 2015. [Online]. Available:
  \url{https://www.medra.org/servlet/aliasResolver?alias=iospress{\&}doi=10.3233/JAE-140148}
\BIBentrySTDinterwordspacing

\bibitem{Darleux2020}
\BIBentryALTinterwordspacing
R.~Darleux, ``{Development of analogous piezoelectric networks for the
  vibration damping of complex structures},'' Ph.D. dissertation, HESAM
  Universit{\'{e}}, 2020. [Online]. Available:
  \url{https://tel.archives-ouvertes.fr/tel-02902755}
\BIBentrySTDinterwordspacing

\bibitem{Giorgio2009}
\BIBentryALTinterwordspacing
I.~Giorgio, A.~Culla, and D.~{Del Vescovo}, ``{Multimode vibration control
  using several piezoelectric transducers shunted with a multiterminal
  network},'' \emph{Archive of Applied Mechanics}, vol.~79, no.~9, pp.
  859--879, sep 2009. [Online]. Available:
  \url{http://link.springer.com/10.1007/s00419-008-0258-x}
\BIBentrySTDinterwordspacing

\bibitem{Yamada2010}
\BIBentryALTinterwordspacing
K.~Yamada, H.~Matsuhisa, H.~Utsuno, and K.~Sawada, ``{Optimum tuning of series
  and parallel LR circuits for passive vibration suppression using
  piezoelectric elements},'' \emph{Journal of Sound and Vibration}, vol. 329,
  no.~24, pp. 5036--5057, nov 2010. [Online]. Available:
  \url{http://dx.doi.org/10.1016/j.jsv.2010.06.021
  https://linkinghub.elsevier.com/retrieve/pii/S0022460X10004116}
\BIBentrySTDinterwordspacing

\bibitem{Thomas2009}
\BIBentryALTinterwordspacing
O.~Thomas, J.-F. De{\"{u}}, and J.~Ducarne, ``{Vibrations of an elastic
  structure with shunted piezoelectric patches: efficient finite element
  formulation and electromechanical coupling coefficients},''
  \emph{International Journal for Numerical Methods in Engineering}, vol.~80,
  no.~2, pp. 235--268, oct 2009. [Online]. Available:
  \url{http://doi.wiley.com/10.1002/nme.2632}
\BIBentrySTDinterwordspacing

\bibitem{Geradin2014}
M.~G{\'{e}}radin and D.~J. Rixen, \emph{{Mechanical vibrations: theory and
  application to structural dynamics}}.\hskip 1em plus 0.5em minus 0.4em\relax
  John Wiley {\&} Sons, 2014.

\bibitem{Grainger1994}
J.~J. Grainger and W.~D.~J. Stevenson, \emph{{Power System Analysis}}.\hskip
  1em plus 0.5em minus 0.4em\relax McGraw-Hill Education, 1994.

\bibitem{Gannett1978}
\BIBentryALTinterwordspacing
J.~Gannett and L.~Chua, ``{Frequency Domain Passivity Conditions for Linear
  Time-Invariant Lumped Networks},'' EECS Department, University of California,
  Berkeley, Berkeley, Tech. Rep. Technical Report No. UCB/ERL M78/21, 1978.
  [Online]. Available:
  \url{https://www2.eecs.berkeley.edu/Pubs/TechRpts/1978/28925.html}
\BIBentrySTDinterwordspacing

\bibitem{Chervov2009}
\BIBentryALTinterwordspacing
A.~Chervov, G.~Falqui, and V.~Rubtsov, ``{Algebraic properties of Manin
  matrices 1},'' \emph{Advances in Applied Mathematics}, vol.~43, no.~3, pp.
  239--315, sep 2009. [Online]. Available:
  \url{https://linkinghub.elsevier.com/retrieve/pii/S0196885809000608}
\BIBentrySTDinterwordspacing

\bibitem{Piefort2001}
\BIBentryALTinterwordspacing
V.~Piefort, ``{Finite Element Modelling of Piezoelectric Active Structures},''
  Ph.D. dissertation, Universit{\'{e}} Libre de Bruxelles, 2001. [Online].
  Available: \url{https://scmero.ulb.ac.be/Publications/Thesis/Piefort01.pdf}
\BIBentrySTDinterwordspacing

\bibitem{Raze2020c}
\BIBentryALTinterwordspacing
G.~Raze, A.~Jadoul, S.~Guichaux, V.~Broun, and G.~Kerschen, ``{A digital
  nonlinear piezoelectric tuned vibration absorber},'' \emph{Smart Materials
  and Structures}, vol.~29, no.~1, p. 015007, jan 2020. [Online]. Available:
  \url{https://iopscience.iop.org/article/10.1088/1361-665X/ab5176}
\BIBentrySTDinterwordspacing

\end{thebibliography}
	\bibliographystyle{IEEEtran}

\end{document}